\documentclass[12pt,leqno]{article}
\usepackage{amsfonts}
\usepackage{amsmath, amsthm, amsfonts, amssymb, color}
\usepackage{mathrsfs}
\usepackage{color}
\newtheorem{thm}{Theorem}[section]

\newtheorem{rem}[thm]{Remark}
\theoremstyle{definition}

\newcommand{\scr}[1]{\mathscr #1}
\definecolor{wco}{rgb}{0.5,0.2,0.3}

\numberwithin{equation}{section} \theoremstyle{remark}

\newcommand{\ua}{\uparrow}

\title{{\bf Distribution Dependent Stochastic Differential Equations}\footnote{Supported in
 part by  NNSFC (11771326, 11831014, 11921001, 11801406).} }
\author{
{\bf  Xing Huang$^{a)}$, Panpan Ren$^{b)}$,    Feng-Yu Wang$^{a)}$  }\\
\footnotesize{$^{a)}$ Center for Applied Mathematics, Tianjin University, Tianjin 300072, China}\\
 \footnotesize{$^{b)}$ Department of Mathematics,
City University of HongKong, HongKong, China}\\
\footnotesize{xinghuang@tju.edu.cn, rppzoe@gmail.com,  wangfy@tju.edu.cn}}
\begin{document}
\allowdisplaybreaks
\def\R{\mathbb R}  \def\ff{\frac} \def\ss{\sqrt} \def\B{\mathbf
B}
\def\N{\mathbb N} \def\kk{\kappa} \def\m{{\bf m}}
\def\ee{\varepsilon}\def\ddd{D^*}
\def\dd{\delta} \def\DD{\Delta} \def\vv{\varepsilon} \def\rr{\rho}
\def\<{\langle} \def\>{\rangle}
  \def\nn{\nabla} \def\pp{\partial} \def\E{\mathbb E}
\def\d{\text{\rm{d}}} \def\bb{\beta} \def\aa{\alpha} \def\D{\scr D}
  \def\si{\sigma} \def\ess{\text{\rm{ess}}}\def\s{{\bf s}}
\def\beg{\begin} \def\beq{\begin{equation}}  \def\F{\scr F}
\def\Ric{\mathcal Ric} \def\Hess{\text{\rm{Hess}}}
\def\e{\text{\rm{e}}} \def\ua{\underline a} \def\OO{\Omega}  \def\oo{\omega}
 \def\tt{\tilde}\def\[{\lfloor} \def\]{\rfloor}
\def\cut{\text{\rm{cut}}} \def\P{\mathbb P} \def\ifn{I_n(f^{\bigotimes n})}
\def\C{\scr C}      \def\aaa{\mathbf{r}}     \def\r{r}
\def\gap{\text{\rm{gap}}} \def\prr{\pi_{{\bf m},\varrho}}  \def\r{\mathbf r}
\def\Z{\mathbb Z} \def\vrr{\varrho} 
\def\L{\scr L}\def\Tt{\tt} \def\TT{\tt}\def\II{\mathbb I}
\def\i{{\rm in}}\def\Sect{{\rm Sect}}  \def\H{\mathbb H}
\def\M{\mathbb M}\def\Q{\mathbb Q} \def\texto{\text{o}} 
\def\Rank{{\rm Rank}} \def\B{\scr B} \def\i{{\rm i}} \def\HR{\hat{\R}^d}
\def\to{\rightarrow}\def\l{\ell}\def\iint{\int}\def\gg{\gamma}
\def\EE{\scr E} \def\W{\mathbb W}
\def\A{\scr A} \def\Lip{{\rm Lip}}\def\S{\mathbb S}
\def\BB{\scr B}\def\Ent{{\rm Ent}} \def\i{{\rm i}}\def\itparallel{{\it\parallel}}
\def\g{{\mathbf g}}\def\Sect{{\mathcal Sec}}\def\T{\mathcal T}\def\BB{{\bf B}}
\def\f{\mathbf f} \def\g{\mathbf g}\def\BL{{\bf L}}  \def\BG{{\mathbb G}}
\def\Bd{{D^E}} \def\BdP{D^E_\phi} \def\Bdd{{\bf \dd}} \def\Bs{{\bf s}} \def\GA{\scr A}
\def\Bg{{\bf g}}  \def\Bdd{\psi_B} \def\supp{{\rm supp}}\def\div{{\rm div}}
\def\ddiv{{\rm div}}\def\osc{{\bf osc}}\def\1{{\bf 1}}\def\BD{\mathbb D}\def\GG{\Gamma}

\maketitle

\begin{abstract} Due to their intrinsic link with nonlinear Fokker-Planck equations and many other applications, distribution dependent stochastic differential equations (DDSDEs for short) have been intensively investigated.
In this paper we summarize some recent progresses  in the study of DDSDEs, which include  the correspondence of weak solutions and nonlinear Fokker-Planck equations,
the well-posedness, regularity estimates,  exponential ergodicity, long time large deviations, and comparison theorems.   \end{abstract} \noindent
 AMS subject Classification:\  60B05, 60B10.   \\
\noindent
 Keywords:   DDSDE,  nonlinear Fokker-Planck equation,  Bismut formula,  Wasserstein distance, gradient estimate.

 \vskip 2cm

 \section{Introduction}



 To characterize nonlinear PDEs in Vlasov's kinetic theory, Kac \cite{Kac, Kac2} proposed the $``$propagation of chaos" of mean field particle systems,  which stimulated  McKean \cite{McKean}  to study   nonlinear Fokker-Planck equations using stochastic differential equations with  distribution dependent  drifts,  see \cite{SN} for a   theory on mean field particle systems and applications.

 In general, a nonlinear Fokker-Planck equation can be characterized by  the following distribution dependent stochastic differential equations (DDSDEs for short):
\beq\label{DDSDE} \d X_t= b(t,X_t,\L_{X_t})\d t+\si(t,X_t,\L_{X_t})\d W_t,\end{equation}
where $W_t$ is an $m$-dimensional Brownian motion on a complete filtration probability space $(\OO, \{\F_t\}_{t\ge 0},\P)$,   $\L_\xi$ is the distribution (i.e. the law) of a random variable $\xi$,
\beg{align*}&b=(b_i)_{1\le i\le d}: [0,\infty)\times\R^d\times \scr P(\R^d)\to \R^d,\\
&\si=(\si_{ij})_{1\le i\le d, 1\le j\le m}: [0,\infty)\times\R^d\times \scr P(\R^d)\to \R^d\otimes\R^m\end{align*}  are measurable, and   $\scr P(\R^d)$ is the space of probability measures on $\R^d$ equipped with the weak topology. Due to the pioneering work \cite{McKean} of McKean,  the DDSDE \eqref{DDSDE} is also called McKean-Vlasov SDE or mean field SDE.

\beg{defn}  Let $s\ge 0$.\beg{enumerate}
\item[$(1)$] A continuous adapted process  $(X_{s,t})_{t\ge s}$   is called a  solution of \eqref{DDSDE} from time   $s$, if
$$\int_s^t\E\big[|b(r,X_{s,r}, \L_{X_{s,r}})|+  \| \si(r,X_{s,r},\L_{X_{s,r}})\|^2\big]\d r<\infty,\ \ t\ge s,$$ and   $\P$-a.s.
  $$X_{s,t} = X_{s,s} +\int_s^t b(r,X_{s,r}, \L_{X_{s,r}})\d r + \int_s^t \si(r,X_{s,r},\L_{X_{s,r}})\d W_r,\ \ t\ge s.$$
When $s=0$ we simply denote   $X_t=X_{0,t}$.
\item[$(2)$] A couple $(\tt X_{s,t},\tt W_t)_{t\ge s}$   is called a  weak solution of \eqref{DDSDE} from time $s$, if   $\tt W_t$  is the   $m$-dimensional Brownian motion on a complete filtration probability space $(\tt\OO,\{\tt\F_t\}_{t\ge 0}, \tt\P)$  such that   $(\tt X_{s,t})_{t\ge s}$  is a solution of \eqref{DDSDE}  from time   $s$
for   $(\tt W_t, \tt\P)$ replacing   $(W_t,\P)$.    \eqref{DDSDE} is called weakly unique for an initial distribution $\nu\in\scr P(\R^d)$, if all weak solutions with distribution $\nu$ at time $s$ are equal in law.
\item[$(3)$] Let $\hat {\scr P}(\R^d)$ be a subspace of $\scr P(\R^d)$. \eqref{DDSDE} is said to have strong (respectively, weak) well-posedness for initial distributions in $\hat {\scr P}(\R^d)$, if for any  $\F_s$-measurable $X_{s,s}$  with $\L_{X_{s,s}}\in \hat {\scr P}(\R^d)$ (respectively, any initial distribution $\nu\in \hat {\scr P}(\R^d)$ at time $s$), it has a unique strong (respectively, weak) solution.
    We call the  equation well-posed   if it is both  strongly and weakly well-posed.\end{enumerate}
\end{defn}
According to Yamada-Watanabe principle, for classical SDEs the strong well-posedness implies the weak one. But this does not apply to
DDSDEs, see Theorem \ref{YM} below for a modified Yamada-Watanabe principle.

In this paper, we summarize the following recent progress in the study of the DDSDE \eqref{DDSDE}:  the correspondence between the weak solution of \eqref{DDSDE} and the associated nonlinear Fokker-Planck equation (Section 2),   criteria on the well-posedness (i.e. existence and uniqueness of solutions) (Section 3),   regularity of distributions (Section 4),  exponential ergodicity (Section 5),  long time large deviations  (Section 6), and comparison theorems (Section 7). Corresponding  results for general models of path-distribution dependent SDEs/SPDEs  can be found in \cite{BRW20, HRW,  RTW20}.

\section{Weak solution and nonlinear Fokker-Planck equation}

In this part, we first introduce the $``$superposition principle" which provides a correspondence between the weak solution of \eqref{DDSDE} and the solution of the associated nonlinear Fokker-Planck equation on $\scr P(\R^d)$, then present some typical examples.

\subsection{Superposition principle}
Consider the following nonlinear Fokker-Planck equation on $\scr P(\R^d)$:
\beq\label{FPE} \pp_t\mu_t = L_{t,\mu_t}^*\mu_t,
\end{equation}
where for any $(t,\mu)\in [0,\infty)\times\scr P(\R^d)$, the   Kolmogorov operator $L_{t,\mu}$  on $\R^d$ is given by
$$L_{t,\mu}:= \ff 1 2 \sum_{i,j=1}^d (\si\si^*)_{ij} (t,\cdot,\mu) \pp_i\pp_j +\sum_{i=1}^d b_i(t,\cdot,\mu) \pp_i,$$
for $\si^*$ being the transposition of $\si$.

\beg{defn}  For $s\ge 0$, $\mu_\cdot\in C([s,\infty); \scr P(\R^d))$ is called a solution of \eqref{FPE} from time $s$, if
$$\int_s^t \d r\int_{\R^d} \big\{\|\si(r,x,\mu_r)\|^2+|b(r,x,\mu_r)|\big\}\mu_r(\d x) <\infty, \  \ t>s,$$
and for any $f\in C_0^\infty(\R^d)$,
\beq\label{LLK} \mu_t(f):=\int_{\R^d} f\d\mu_t= \mu_s(f)+\int_s^t \mu_r(L_{r,\mu_r}f)\d r,\ \ t\ge s.\end{equation}
\end{defn}

Now, assume that $(\tt X_t,\tt W_t)_{t\ge s}$ is a weak solution of \eqref{DDSDE} from time $s$ under a complete filtration probability space $(\tt\OO, \{\tt \F_t\}_{t\ge s}, \tt\P)$, and let
$\mu_t= \L_{\tt X_t|\tt\P}:= \tt\P\circ (\tt X_t)^{-1}$ be the   distribution of $\tt X_t$ under the probability $\tt\P$. By It\^o's formula we have
$$ \d f(\tt X_t)=\big\{ L_{t,\mu_t} f(\tt X_t) \big\}\d t +\<\nn f(\tt X_t), \si(t,\tt X_t, \mu_t)\d \tt W_t\>.$$
Integrating both sides  over $[s,t]$ and  taking expectations, we obtain \eqref{LLK} so that $\mu_\cdot$ solves \eqref{FPE} by definition.
Indeed,   the following result due to \cite{RWBR18,RWBR} also ensures the converse, i.e. a solution of \eqref{FPE}  gives a weak solution of \eqref{DDSDE},
see  Section 2 of \cite{RWBR} (and \cite{RWBR18}).

\begin{thm} [\cite{RWBR18,RWBR}]
Let $(s,\zeta)\in[0,\infty)\times\scr P(\R^d)$. Then the DDSDE \eqref{DDSDE} has a weak solution $(\tt X_t,\tt W_t)_{t\ge s}$  starting from $s$ with   $\L_{\tt X_s|\tt\P}=\zeta$, if and only if \eqref{FPE} has a solution $(\mu_t)_{t\geq s}$ starting from $s$ with $\mu_s=\zeta$. In this case $\mu_t=\scr L_{\tt X_t| \tt\P}$, $t\geq s$.
\end{thm}

\subsection{Some examples}
In this part, we introduce some  typical nonlinear PDES and state their corresponding DDSDEs.
\paragraph{Example 2.1 (Landau type equations).}Consider the following nonlinear PDE  for probability density functions $(f_t)_{t\ge 0}$  on $\R^d$:
 \beq\label{LD} \pp_t f_t=\ff 1 2 \text{div}\bigg\{\int_{\R^d} a(\cdot-z)\big(f_t(z)\nn f_t-f_t\nn f_t(z)\big)\d z\bigg\},\end{equation}
where   $a:\R^d\to \R^d\otimes\R^d$ has weak derivatives.
For the real-world model of homogenous Landau equation we have  $d=3$ and
  $$a(x)= |x|^{2+\gg} \Big(I-\ff{x\otimes x}{|x|^2}\Big),\ \ x\in \R^3 $$  for some constant   $\gg\in [-3,1].$
In this case \eqref{LD}  is a limit version  of Boltzmann equation (for thermodynamic system) when all collisions become grazing.
To characterize this equation using SDE,  let   $m=d$, $b= \ff 1 2\text{div} a$ and  $ \si= \ss a$. Consider the  DDSDE
  \beq\label{SD0} \d X_t= (b*\L_{X_t})(X_t)\d t +(\si *\L_{X_t})(X_t)\d B_t,\end{equation}   where
  $$(f*\mu)(x):=\int_{\R^d} f(x-z)\mu(\d z).$$  Then the distribution density
  $ f_t(x):= \ff{ \L_{X_t}(\d x)}{\d x} $   solves the Landau type  equation   \eqref{LD}.

\paragraph{Example 2.2 (Porous media equation).}Consider the following  nonlinear PDE for probability density functions on $\R^d$:
\beq\label{PM}   \pp_t f_t= \DD f_t^3.\end{equation}
Then  for any solution to the DDSDE \eqref{DDSDE} with coefficients
$$b=0,\ \ \si(x,\mu)= \ss 2 \ff{\d\mu}{\d x}(x) I_{d\times d},$$
  the probability density function  solves the porous media equation \eqref{PM}.

\paragraph{Example 2.3 (Granular media equation).}Consider the following nonlinear PDE for probability density functions on $\R^d$:
\beq\label{GM}   \pp_t f_t= \DD f_t +{\rm div} \big\{f_t\nn V + f_t \nn (W*f_t)\big\}. \end{equation}
Then the associated  DDSDE \eqref{DDSDE} has coefficients
 $$b(x,\mu)= -\nn V(x) -\nn (W*\mu)(x) ,\ \ \si(x,\mu)= \ss 2 I_{d\times d},$$
where
  $$(W*\mu)(x):=\int_{\R^d} W(x-y)\mu(\d y).$$

\section{Well-posedness }
We first introduce  a fixed-point argument in distribution and a modified Yamada-Watanabe principle, then present results on the existence and uniqueness for monotone and singular coefficients respectively.

\subsection{Fixed-point in distribution and  Yamada-Watanabe principle}
Let $\hat {\scr P}(\R^d)$ be a subspace of $\scr P(\R^d)$, and let $\hat\rr$ be a complete metric on $\hat {\scr P}(\R^d)$ inducing the Borel sigma algebra of the weak topology.  Typical examples include
$$\hat{\scr P}(\R^d)= \scr P_p(\R^d):=\big\{\mu\in \scr P(\R^d): \mu(|\cdot|^p)<\infty\big\}$$ for $p>0$, with $L^p$-Wasserstein distance
$$\W_p(\mu,\nu):= \inf_{\pi\in \C(\mu,\nu)} \bigg(\int_{\R^d\times\R^d} |x-y|^p \pi(\d x,\d y) \bigg)^{\ff 1{p\lor 1}},\ \ \mu,\nu\in \scr P_p(\R^d).$$ When $p=0$ this reduces to the total variation norm
$$\|\mu-\nu\|_{TV}:=2\sup_{A\in\B(\R^d)} |\mu(A)-\nu(A)|.$$
For any $T>s\ge 0$ and $\nu\in \hat{\scr P}(\R^d)$, consider the path space over $\hat{\scr P}(\R^d)$
$$\hat \C_{s,T}^\nu:=\big\{\mu_\cdot\in C([s,T]; \hat{\scr P}(\R^d)): \mu_s=\nu\big\},$$
which is then complete under the metric
$$\hat\rr_{s,T}(\mu_\cdot,\nu_\cdot):= \sup_{t\in [s,T]} \hat\rr(\mu_t,\nu_t).$$

\beg{thm} \label{TFIX}  Let $T>s\ge 0$,  and let $X_s$ be an $\F_s$-measurable random variable with $\nu:=\L_{X_s}\in \hat{\scr P}(\R^d).$   Assume that for any
$\mu\in \hat \C_{s,T}^\nu,$ the classical SDE
\beq\label{CSDE} \d X_t^\mu=b(t,X_t^\mu, \mu_t)\d t +\si(t, X_t^\mu, \mu_t)\d W_t,\ \ t\in [s,T], X_s^{\mu}=X_s\end{equation}
has a unique solution, and  the map
$$\mu\in \hat \C_{s,T}^\nu \mapsto \Phi_{s,T} \mu := (\L_{X_t^\mu})_{t\in [s,T]}\in \hat{\C}_{s,T}^\nu$$ is   contractive. Then the DDSDE $\eqref{DDSDE}$ has well-posedness for initial distributions in $\hat{\scr P}(\R^d).$
\end{thm}

\beg{proof} By the fixed-point theorem, the map $\Phi_{s,T}$ has a unique fixed point $\mu$ in $\hat \C_{s,T}^\mu$, so that by the definition of $\Phi_{s,T}$ we have
$\L_{X_t^\mu}=\mu_t, t\in [s,T]$, i.e. in this case $(X_t^\mu)_{t\in [s,T]}$ is a solution of \eqref{DDSDE} from time $s$ starting at $X_s.$ If \eqref{DDSDE} has another solution $(\hat X_t)_{t\in [s,T]}$ with $\L_{\hat X_\cdot}\in \hat C_{s,T}^\mu$, then $\mu:=\L_{\hat X_\cdot}$ is a fixed point  of $\Phi_{s,T}$ so that $\L_{X_t}=\L_{\hat X_t}=:\mu_t, t\in [s,T].$
Therefore, by the uniqueness of \eqref{CSDE} we have $\L_{X_t}=\L_{\hat X_t}=X_t^\mu$, which implies the uniqueness of \eqref{DDSDE} with $\L_{X_\cdot}\in \hat \C_{s,T}^\mu$.
Since the strong well-posedness of \eqref{CSDE} implies the weak one, the same argument leads to the weak well-posedness of the DDSDE \eqref{DDSDE} starting from $\nu$ at time $s$.
\end{proof}

The Yamada-Watanabe principle \cite{YW} (see \cite{Kurt} for a general version) is a fundamental tool in the study of well-posedness for SDEs with singular coefficients. In the present distribution dependent setting, the original statement does not apply, but we have the following modified version due to \cite[Lemma 3.4]{HW19}.

\beg{thm}[\cite{HW19}] \label{YM}   Let $T>s\ge 0$,  and let $X_s$ be an $\F_s$-measurable random variable with $\nu:=\L_{X_s}\in \hat{\scr P}(\R^d).$   Assume that for any
$\mu\in \hat \C_{s,T}^\nu,$ the classical SDE $\eqref{CSDE}$ has a unique solution with  initial value $X_s$ at time $s$. If \eqref{DDSDE} for $t\in [s,T]$  has a weak solution with initial distribution $\nu$ at time $s$, and has pathwise uniqueness with initial value $X_s$ at times $s$, then it has well-posedness for initial distributions in  $\hat{\scr P}(\R^d).$  \end{thm}

\subsection{The monotone case}
\beg{enumerate} \item[$(H_3^1)$]   For every $t\ge 0$, $b_t$ is continuous on $\R^d\times\scr P_\theta(\R^d)$, $b$ is bounded on bounded sets in $[0,\infty)\times \R^d\times \scr P_\theta(\R^d)$. Moreover, there exists
$K \in L^1_{loc}([0,\infty);(0,\infty))$ such that
 \beg{align*}  &\|\si(t,x,\mu)-\si(t,y,\nu)\|^2\le K(t)\big\{|x-y|^2+\W_\theta(\mu,\nu)^2\big\},\\
 &\<b(t,x,\mu)-b(t,y,\nu),x-y\>
   \le K(t)\big\{|x-y|^2 + \W_\theta(\mu,\nu) |x-y|\big\},\\
   &|b(t,0,\dd_0)|+\|\si(t,0,\dd_0)\|_{HS}^2 \le K(t),\  \ t\ge 0, x,y\in \R^d, \mu,\nu\in
\scr P_\theta(\R^d),\end{align*} where $\dd_0$ is the Dirac measure at $0\in\R^d.$
\end{enumerate}
Under this monotone condition we have the following result essentially due to \cite{W18}, where a stronger growth condition
on $|b(t,0,\mu)|$ is assumed.
See also \cite{HSS19} for the well-posedness under integrated Lyapunov
conditions which may cover more examples.

\beg{thm}[\cite{W18}]\label{T1.1dd} Assume $(H_3^1)$  for some $\theta\in [1,\infty)$, and let $\si(t,x,\mu)$ does not depend on $\mu$    when $\theta<2$.
 \beg{enumerate} \item[$(1)$] The DDSDE $\eqref{DDSDE}$ has well-posedness for initial distributions in $\scr P_\theta(\R^d)$.
   Moreover, for any $p\ge \theta$ and $s\ge 0$, $\E|X_{s,s}|^{p}<\infty$ implies
\begin{equation*}\E \sup_{t\in [s,T]}|X_{s,t} |^{p}<\infty,\ \ T\ge t\ge s\ge  0.\end{equation*}
\item[$(2)$] There exists increasing $\psi: [0,\infty)\to [0,\infty)$ such that for any two solutions $X_{s,t}$ and $Y_{s,t}$ of $\eqref{DDSDE}$ with $\L_{X_{s,s}},\L_{Y_{s,s}}\in \scr P_\theta(\R^d)$,
\beq\label{CCdd} \E |X_{s,t}-Y_{s,t}|^\theta\le \big(\E|X_{s,s}-Y_{s,s}|^\theta\big)\e^{ \int_s^t\psi(r)\d r},\ \ t\ge s\ge 0.\end{equation} Consequently,
\beq\label{CC2}\lim_{\E|X_{s,s}-Y_{s,s}|^\theta\to 0} \P\Big(\sup_{r\in [s,t]}|X_{s,r}-Y_{s,r}|\ge \vv\Big)=0,\ \ t>s\ge 0,\vv>0;\end{equation} and
\beq\label{X0}  \W_\theta(P_{s,t}^*\mu_0, P_{s,t}^*\nu_0)^\theta\le   \W_2(\mu_0,\nu_0)^\theta\e^{\int_s^t\psi(r)\d r},\ \ t\ge s\ge 0.\end{equation}
\end{enumerate} \end{thm}

\beg{proof} We briefly explain the proof of Theorem \ref{T1.1dd}(1), while (2) can be easily proven by using It\^o's formula. For any $T>s\ge 0$, $\nu\in \scr P_\theta(\R^d)$ and $\mu\in \hat\C_{s,T}^\nu$,  $(H_3^1)$ implies that \eqref{CSDE} for $t\in [s,T]$ is well-posed with initial distribution $\nu$ at $s$.
Moreover, by It\^o's formula, and $(H_3^1)$ with $\si(t,x,\mu)$ not depending on $\mu$ when $\theta<2$, we find a large enough constant
$\lambda>0$ such that $\Phi_{s,T}$ is contractive on $\hat \C_{s,T}^\nu$ under the complete metric
$$\hat\rr_{s,T}(\mu,\tt\mu):= \sup_{t\in [s,T]} \e^{-\lambda (t-s)} \W_\theta(\mu_t,\tt\mu_t),\ \ \mu,\tt\mu\in \hat C_{s,T}^\nu.$$
Then the well-posedness follows from Theorem \ref{TFIX}. \end{proof}

\subsection{The singular case}

In this part, we consider the existence and uniqueness of  \eqref{DDSDE} with singular drift and non-degenerate noise.
We first introduce  some results derived in  \cite{Zhao, RZ, Z4} for distribution dependent drifts satisfying local integrability conditions in time and space but bounded in distribution,
  in  \cite{HW20c}   for  the case with locally integrable drifts having linear growth in distribution, and in \cite{HW20a} for drifts with
an integrable term and a Lipchitz term. These three situations are mutually incomparable.

\subsubsection{Integrability  in time-space and boundedness in distribution}
When the noise is possibly degenerate, the strong/weak well-posedness will be discussed in the next section under a monotone condition.

We will consider weak solutions  having finite  $\phi$-moment,  for $\phi$ in the following class:
$$  \mathbf \Phi:=\big\{\phi\in C^\infty([0,\infty); [1,\infty)): 0\le \phi'\le c\phi\ \text{for\ some\ constant\ } c>0\big\}.$$
Let
$$\scr P_\phi(\R^d):=\big\{\mu\in \scr P(\R^d): \|\mu\|_\phi:= \mu(\phi(|\cdot|))<\infty\big\},$$
which is equipped with the $\phi$-total variation norm
$$\|\mu-\nu\|_{\phi,TV} := \sup_{|f|\le \phi(|\cdot|)} \big|\mu(f)-\nu(f)\big|,\ \ \mu,\nu\in \scr P(\R^d).$$
We denote $\|\cdot\|_{\phi,TV}$ by $\|\cdot\|_{\theta,TV}$ when $\phi=1+|\cdot|^\theta$ for some $\theta\geq 0$.
For fixed $T>0$, let
$$\C_{T,\phi}= C([0,T]; \scr P_\phi(\R^d)):=\Big\{\mu: [0,T]\to \scr P_\phi(\R^d), \lim_{t\to s} \|\mu_t-\mu_s\|_{\phi,TV}=0,\ s\in [0,T]\Big\},$$
which is a complete  space under  the metric
$$\rr_{T,\phi}(\mu,\nu):= \sup_{t\in [0,T]} \|\mu_t-\nu_t\|_{\phi,TV}.$$
For any $\mu\in\scr C_{T,\phi}$, denote
$$b^\mu(t,x):= b(t,x,\mu_t),\ \ \si^\mu(t,x):= \si(t,x,\mu_t), \ \ a^\mu(t,x):= \ff 1 2 \big\{\si^\mu(\si^\mu)^*\big\}(t,x),\ \ (t,x)\in [0,T]\times\R^d.$$
\beg{defn}[Linear Functional Derivative] Let $\phi\in \mathbf \Phi.$ A function $f: \scr P_\phi(\R^d)\to\R$
is said to have  linear functional derivative
$D^Ff: \scr P_\phi(\R^d)\to \R,$ if it is measurable, and \beg{enumerate} \item[(i)] $D^Ff$ is measurable with
$\int_{\R^d}  D^Ff(\mu)\d\mu=0;$
\item[(ii)] For any compact  $K\subset \scr P_\phi(\R^d)$,
$\sup_{\mu\in K} |D^Ff(\mu)|\le k\phi(|\cdot|)$ holds for some constant $k>0$;
\item[(iii)] For any $\mu,\nu\in \scr P_\phi(\R^d),$
$$\lim_{s\downarrow 0} \ff{f((1-s)\mu+s\nu)-f(\mu)}s =\int_{\R^d} D^Ff(\mu)(y)(\nu-\mu)(\d y).$$\end{enumerate} \end{defn}

By taking $\nu=\dd_y$, we see that if $f$ has linear functional derivative, then the convex extrinsic derivative
$$\tt D^E f(\mu)(y):=\lim_{s\downarrow 0} \ff{f((1-s)\mu+s\dd_y)-f(\mu)}s =
D^Ff(\mu)(y)-\int_{\R^d} D^Ff(\mu)(y)\d\mu$$
exists. See \cite{RW19} for links of more derivatives in measure. For $i=1,2$, let
\begin{align}\label{Ti}
\scr I_i=\Big\{(p,q)\in(1,\infty)\times (1,\infty):\frac{d}{p}+\frac{2}{q}<i\Big\}.
\end{align}
\beg{defn}  For any $p\ge 1$,  let $\tt L_p$ be the space of all measurable functions $g$ on $\R^d$ such that
$$\|g\|_{\tt L_p}:=\sup_{z\in \R^d}\bigg(\int_{\R^d} |g(x)|^p1_{\{|x-z|\le 1\}} \d x\bigg)^{\ff 1 p}<\infty.$$
Moreover, for any $p,q\ge 1$, let $\tt L_p^q(T)$ be the space of measurable functions $f$ on $[0,T]\times\R^d$ such that
 $$\|f\|_{\tt L_p^q(T)}:=\sup_{z\in\R^d} \bigg(\int_0^T  \left(\int_{\R^d} |f(t,x)|^p1_{\{|x-z|\le 1\}} \d x\right)^\frac{q}{p}  \d t \bigg)^{\ff 1q}  <\infty.$$
   \end{defn}
It is clear that $$\|f\|_{\tt L_p^q(T)}\leq \|f\|_{L^q([0,T];\tilde{L}_p)}:= \left(\int_0^T\|f(t,\cdot)\|^q_{\tilde{L}_p}\right)^\frac{1}{q}.$$ The following result  is due to \cite[Theorems 3.5 and 3.9]{Zhao}, see also \cite{Sch} for a special case where $\phi(r)=r^2$ and $b(t,x,\cdot)$ is bounded and Lipschitz continuous in the total variation norm uniformly in $(t,x)$.

\beg{thm}[\cite{Zhao}] \label{Zhao} Let $\si\si^*$ be invertible, $\phi\in \mathbf \Phi$, and  $p,q\in (1, \infty)$ with $\vv:= 1-\ff d p -\ff 2 q>0.$
\beg{enumerate} \item[$(1)$] If there exist constants $\aa\in (0,1), N>1, $ and $r>\ff 2 \vv$ such that for any $\mu\in \C_{T,\phi}$,
 $$\sup_{t\in[0,T], x\neq y}\frac{\|a^\mu(t,x)-a^\mu(t,y)\|}{|x-y|^\alpha}+\sup_{(t,x)\in [0,T]\times\R^d} \big\{\|a^\mu\|+\|(a^\mu)^{-1}\|\big\}(t,x) +\|b^\mu\|_{\tt L_p^q(T)}  \le N,$$
and that
$$ \lim_{\rr_{T,\phi}(\nu,\mu)\to 0}  \bigg\{ \int_0^T \|a^\mu(t,\cdot)-a^\nu(t,\cdot)\|_\infty^r \d t  +\|b^\mu-b^\nu\|_{\tt L_p^q(T)}   \bigg\}=0,  $$
then  $\eqref{DDSDE}$ has a weak solution for $t\in [0,T]$ and   any initial distribution in $\scr P_\phi(\R^d)$.
\item[$(2)$]  In addition to conditions in $(1)$,
if  for any $(t,x)\in [0,T]\times\R^d$, $\si(t,x,\cdot)$ has linear functional derivative on $\scr P_\phi(\R^d)$, and there exist constants $\bb\in(0,1), C>0$ and some $K\in L^q([0,T]; (0,\infty))$ such that
\beg{align*} & \sup_{(t,\mu)\in [0,T]\times \scr P_\phi(\R^d)} \big\|D^F\si(t,x,\cdot)(\mu)(y)-D^F\si(t,x',\cdot)(\mu)(y')\big\|\\
&\le C(|x-x'|+|y-y'|)^\bb,
 \qquad\  x,x',y,y'\in\R^d,\\
& \|b(t,\cdot,\mu)-b(t,\cdot,\nu)\|_{\tt L_p} \le K(t) \|\mu-\nu\|_{\phi,TV} ,\ \ t\in [0,T], \mu,\nu\in \scr P_\phi(\R^d),\end{align*}
then   $\eqref{DDSDE}$ is has weak well-posedness for $t\in [0,T]$ and initial
 distribution in  $\scr P_\phi(\R^d)$.\end{enumerate}
\end{thm}
When $\sigma=\sqrt{2}I_{d\times d}$ and
 $$|b(t,x,\mu)|\le \int_{\R^d} h_t(x-y)\mu(\d y)$$ holds for some $(p,q)\in \scr I_1$ and  $h\ge 0$ with $\|h\|_{L^q([0,T]; \tt L_p)}<\infty,$     the   well-posedness for \eqref{DDSDE} is proved in \cite[Theorem 1.1]{RZ}.
In general,    \cite{RZ}  presents the following result.

\begin{thm} [\cite{RZ}] \label{RZ}   Assume that for each $t,x$, $b(t,x,\cdot)$ and $\sigma(t,x,\cdot)$ are weakly continuous, and there exist
$c_0 > 1$ and $\gamma\in(0, 1]$ such that for all $t \geq  0, x, y,\xi\in\R^d$ and $\mu\in \scr P(\R^d)$,
$$c_0^{-1}|\xi|\leq |\sigma(t,x,\mu)\xi|\leq c_0|\xi|, \ \ |\sigma(t,x,\mu)-\sigma(t,y,\mu)|\leq |x-y|^\gamma.$$ Moreover, under the weak topology of $\scr P(\R^d)$,
$$\sup_{\mu\in C([0,T]; \scr P(\R^d))}\|b^\mu\|_{\tilde{L}_p^q(T)}<\infty$$
holds for some $(p,q)\in \scr I_1$.
Then for any $\beta > 2$ and $ \nu\in\scr P_\beta(\R^d)$, there exists a weak solution to \eqref{DDSDE} with initial distribution $\nu$.
If in addition, $\sigma(t,x,\mu)$ does not depend on $\mu$, $|\nabla\sigma|\in \tilde{L}^{q_1}_{p_1}(T)$
and
$$\|b(t,\cdot,\mu)-b(t,\cdot,\nu)\|_{\tilde{L}_p}\leq \ell_t\|\mu-\nu\|_{\theta,TV},\ \ \mu,\nu\in\scr P_\theta$$
for some $\ell\in L^q([0,T])$, $\theta\geq 1$ and $(p_1,q_1)\in\scr I_1,$
then for any $\beta>2\theta$,  \eqref{DDSDE} has well-posedness from time $0$  for   initial distributions in $\scr P_\bb(\R^d)$.
\end{thm}
The following weak existence for \eqref{DDSDE} with   supercritical drift is due to \cite{Z4}.
\begin{thm} [\cite{Z4}] \label{Z4} Let $\sigma=\sqrt{2}I_{d\times d}$,
  $b(t,x,\mu)=\int_{\R^d}K(t,x,y)\mu(\d y)$ for some measurable function $K$ on $[0,T]\times\R^d\times\R^d$ such that  $\mathrm{div }K(t,\cdot,y)\leq 0$ and
$$K(t,x,y)\leq h_t(x,y)$$ holds for    some $(p,q)\in\scr I_2$ and $h\ge 0$ with $\|h\|_{L^q([0,T]; \tt L_p)}<\infty.$
Then for any $\beta\in[0,2/(\frac{d}{p}+\frac{2}{q}))$ and $\nu\in \scr P_{\beta}(\R^d)$, \eqref{DDSDE} has a weak solution with initial distribution $\nu$.
\end{thm}

\subsubsection{ Integrability in time-space with linear growth in distribution }
Comparing with above results, besides the singularity in $x$ in the following we  also allow $b(t,x, \mu)$ to have a linear growth in $\mu$.

\beg{enumerate} \item[$(H_3^2)$]   Let $\theta\ge  1.$
    There exists a constant $K>0$ such that  for any $t\in [0,T], x,y\in\R^d$ and $\mu,\nu \in \scr P_\theta$,
\beg{align*} &\|\si(t,x,\mu)\|^2\lor \|(\si\si^*)^{-1}(t,x,\mu)\|\le K,\\
&\|\si(t,x,\mu)-\si(t,y,\nu)\|\le K\big(|x-y|+\W_\theta(\mu, \nu)\big),\\
&\|\{\si(t,x,\mu)-\si(t,y,\mu)\}-\{\si(t,x,\nu)-\si(t,y,\nu)\}\|\le K|x-y|\W_\theta(\mu, \nu).
\end{align*}
Moreover,   there exists nonnegative
$f\in \tilde{L}_p^q(T)$ for some $(p,q)\in \scr I_1$
such that
\beg{align*}  &|b(t,x,\mu)|\le  (1+\|\mu\|_\theta)f_t(x),\\
&|b(t,x,\mu)-b(t,x,\nu)|\le  f_t(x)\|\mu-\nu\|_{\theta,TV},\ \ t\in [0,T], x\in\R^d, \mu,\nu\in \scr P_\theta.\end{align*}
\end{enumerate}
 \beg{thm}[\cite{HW20c}]\label{T1dd} Assume $(H_3^2)$. Then $\eqref{DDSDE}$ is well-posed for initial distributions in $\scr P_{\theta+}:=\cap_{m>\theta} \scr P_m,$ and  the solution satisfies $\L_{X_\cdot}\in C([0,T];\scr P_\theta),$
 the space of continuous maps from $[0,T]$ to $\scr P_\theta$ under the metric $\W_\theta$. Moreover,
\beq\label{MM} \E\Big[ \sup_{t\in [0,T]} |X_t|^\theta \Big]<\infty.\end{equation}
\end{thm}

\subsubsection{Drifts with  time-space integrable and    Lipschitz terms}
In this part we allow the drift to include a Lipschitz continuous term in $x$, but the price we have to pay is that the singular term is in $L_p^q(T)$ rather than $\tt L_p^q(T)$  and the diffusion does not depend on distribution.

 For any $p,q\ge 1$, let $L_p^q(T)$ be the space of measurable functions $f$ on $[0,T]\times\R^d$ such that
 $$\|f\|_{L_p^q(T)}:= \bigg(\int_0^T  \left(\int_{\R^d} |f(t,x)|^p \d x\right)^\frac{q}{p}  \d t \bigg)^{\ff 1q}  <\infty.$$
\begin{enumerate}
\item[$(H_3^{3a})$] $\si(t,x,\mu)=\si(t,x)$ does not depend on $\mu$ and is uniformly continuous in $x\in\R^d$ uniformly in $t\in [0,T];$ the weak gradient $\nn \si(t,\cdot)$ exists for a.e. $t\in [0,T]$ satisfying
    $|\nn \si|^2\in L_p^q(T)$ for some  $(p,q)\in \scr I_1$;  and there exists a constant $K_1\ge 1$ such that
\begin{align} \label{si-} K_1^{-1} I_{d\times d} \le (\sigma\sigma^\ast)(t,x)\le K_1 I_{d\times d},\ \ (t,x)\in [0,T]\times\R^d.\end{align}
 \item[$(H_3^{3b})$]
$b=\bar{b}+\hat{b}$, where $\bar b$ and $\hat b$ satisfy
\begin{equation}\label{con}\beg{split}
&|\hat{b}(t,x,\gamma)-\hat{b}(t,y,\tilde{\gamma})|+|\bar{b}(t,x,\gamma)-\bar{b}(t,x,\tilde{\gamma})|\\
&\le K_2(\|\gamma-\tilde{\gamma}\|_{TV}+\W_\theta(\gamma,\tilde{\gamma})+|x-y|),\ \ t\in[0,T], x,y \in \R^d, \gamma,\tilde{\gamma}\in \scr P_\theta(\R^d)\end{split}
\end{equation}for some constants $\theta,K_2\ge 1,$ and
for $(p,q)$ in $(H_3^{3a})$, it holds that
\beq\label{ihp} \sup_{t\in [0,T],\gamma\in \scr P_{\theta}(\R^d)} |\hat{b}(t,0,\gamma)|+\sup_{\mu\in C([0,T]; \scr P_{\theta}(\R^d))}\| |\bar{b}^\mu|^2\|_{ L_{p}^{q}(T)}<\infty.\end{equation}
\item[$(H_3^{3c})$] For any $\mu\in \B([0,T];\scr P(\R^d))$, the class of measurable maps from $[0,T]$ to $\scr P(\R^d)$, $|b^\mu|^2\in L_{p,loc}^q(T)$ for $(p,q)$ in $(H_3^{3a})$. Moreover, there exists an increasing  function $\Gamma:[0,\infty)\to (0,\infty)$ satisfying $\int_{1}^\infty\frac{1}{\Gamma(x)}\d x=\infty$ such that
\begin{align}\label{ne}\<b(t,x,\delta_0),x\>\leq \Gamma(|x|^2),\ \ t\in[0,T],x\in \R^d.\end{align}
In addition, there exists a constant $K_3\ge 1$ such that
 \begin{equation}\label{cono}\beg{split}
&|b(t,x,\gamma)-b(t,x,\tilde{\gamma})|\le K_3\|\gamma-\tilde{\gamma}\|_{TV},\ \ t\in[0,T], x\in \R^d, \gamma,\tilde{\gamma}\in \scr P(\R^d).
\end{split}
\end{equation}
\end{enumerate}
\begin{thm}[\cite{HW20a}] \label{T1.1} Assume $(H_3^{3a})$.
\beg{enumerate}
\item[$(1)$] If $(H_3^{3c})$ holds, then  \eqref{DDSDE}  is   well-posed for   initial  initial distributions in $\scr P_\theta(\R^d)$. Moreover,
\beq\label{B1} \|P_t^*\mu_0- P_t^*\nu_0\|_{TV}^2\le 2\e^{\ff{K_1K_3^2t}2}\|\mu_0-\nu_0\|_{TV}^2,\ \ t\in [0,T], \mu_0,\nu_0\in \scr P(\R^d).\end{equation}
 \item[$(2)$]   Let $(H_3^{3b})$ hold. Then $\eqref{DDSDE}$ is   well-posed  for initial  distributions in $\scr P_{\theta}(\R^d)$. Moreover,   for any  $m\in (\ff \theta 2,\infty)\cap [1,\infty)$,   there exists a constant $c>0$ such that
\beq\label{B2}\begin{split} &\|P_t^*\mu_0- P_t^*\nu_0\|_{TV} +\W_\theta(P_t^*\mu_0, P_t^*\nu_0)\\
&\le c \big\{\|\mu_0-\nu_0\|_{TV}+\W_{2m}(\mu_0,\nu_0)\big\}, \ \ t\in [0,T], \mu_0,\nu_0\in \scr P_{\theta}(\R^d).\end{split}\end{equation}
 \end{enumerate}
\end{thm}
\section{Regularity estimates}

In this section, we introduce some     results on the regularity of distributions for  the DDSDE \eqref{DDSDE}. We first establish the log-Harnack inequality, which implies the $``$gradient estimate"  and entropy estimate, then establish the Bismut formula for the Lions derivative of the distribution, and finally study the derivative estimate
on the distribution. In the first two cases   the noise does not depend on the distribution,  while  the last part applies also to distribution dependent noise.

\subsection{Log-Harnack inequality}

The dimension-free  Harnack inequality  was founded  in  \cite{W97} for diffusion semigroups on Riemannian manifolds, and as a  weaker version the log-Harnack inequality was introduced in   \cite{RW10,W10} for (reflecting) diffusion processes and SDEs.  See the monograph \cite{Wbook} for the study of these type inequalities and  applications.
In  this part, we introduce the log-Harnack inequality established in \cite{W18} and \cite{RW20b} for DDSDEs with non-degenerate and degenerate noise respectively. We will only consider
 distribution, independent noise, since the log-Harnack inequality is not yet available for DDSDEs with distribution dependent noise.

\subsubsection{The non-degenerate case}
Consider the following special version of \eqref{DDSDE}:
\beq\label{E11} \d X_t= b(t,X_t,\L_{X_t})\d t +\si(t,X_t)\d W_t,\end{equation}  where $b$ and $\si$ satisfy the following assumption.

\beg{enumerate} \item[$(H_4^1)$] $\si(t,x)$ is invertible and  Lipschitzian in $x$  locally uniformly in $t\ge 0$,   and there exist increasing functions $\kk_0,\kk_1,\kk_2,\lambda: [0,\infty)\to (0,\infty)$ such that for any $t\in [0,T], x,y\in\R^d$ and $\mu,\nu\in \scr P_2(\R^d)$, we have
\beq\label{A1}\|\si(t,\cdot)^{-1}\|_\infty\le \lambda(t),\ \ |b(t,0,\mu)|^2+\|\si(t,x)\|^2\le \kk_0(t)(1+|x|^2+\mu(|\cdot|^2)),\end{equation}
  \beq\beg{split}\label{A2}& 2\<b(t,x,\mu)-b(t,y,\nu), x-y\> +\|\si(t,x)-\si(t,y)\|_{HS}^2 \\
  &\le \kk_1(t)|x-y|^2 +\kk_2(t) |x-y|\W_2(\mu,\nu).\end{split}
 \end{equation}\end{enumerate}
Obviously, $(H_4^1)$ implies  assumptions   $(H_3^1)$ for $\theta=2$, so that  Theorem \ref{T1.1dd} ensures the well-posedness of \eqref{E11} with initial distributions in $\scr P_2(\R^d)$.
For any $f\in \B_b(\R^d)$, consider
$$P_{s,t} f(\mu):= \E^\mu f(X_{s,t}) =\int_{\R^d} f(y)(P_{s,t}^*\mu)(\d y),\ \ \mu\in \scr P_2(\R^d), t\ge s\ge 0,$$
where $\E^\mu$ is the expectation taking for the solution $(X_{s,t})_{t\ge s}$ of \eqref{E11} with $\L_{X_{s,s}}=\mu$, recall that in this case we denote $P_{s,t}^*\mu= \L_{X_{s,t}}$.
Let
$$\phi (s,t) = \lambda(t)^2 \bigg(\ff{\kk_1(t)}{ 1-\e^{-\kk_1(t)(t-s)} }+ \ff{t\kk_2(t)^2\exp[2(t-s)(\kk_1(t)+\kk_2(t))]}{2 }\bigg),\ \ 0\le s<t.$$

\beg{thm}[\cite{W18}]\label{T3.1} Assume $(H_4^1)$ and let $t>s\ge 0$.
Then for any   $\mu_0,\nu_0\in \scr P_2(\R^d)$,
$$(P_{s,t}\log f)(\nu_0)\le  \log (P_{s,t}f)(\mu_0)+ \phi (s,t)\W_2(\mu_0,\nu_0)^2,\ \ f\in \B_b^+(\R^d).$$
Consequently,   the following assertions hold:
\beg{enumerate} \item[$(1)$] For any $\mu_0,\nu_0\in \scr P_2(\R^d)$,
 $$\|P_{s,t}^*\mu_0-P_{s,t}^*\nu_0\|_{TV}
 \le \ss{2\phi (s,t)} \W_2(\mu_0,\nu_0).$$
 \item[$(2)$] For any $\mu_0,\nu_0\in\scr P_2(\R^d)$, $P_{s,t}^*\mu_0$ and $P_{s,t}^*\nu_0$ are equivalent and the Radon-Nykodim derivative satisfies the entropy estimate
$$ {\rm Ent} (P_{s,t}^*\nu_0|P_{s,t}^*\mu_0):=\int_{\R^d} \bigg\{\log \ff{\d P_{s,t}^*\nu_0}{\d P_{s,t}^*\mu_0}\bigg\}\d P_{s,t}^*\nu_0\le \phi (s,t) \W_2(\mu_0,\nu_0)^2.  $$
\end{enumerate}
\end{thm}

 \beg{proof}[Idea of Proof]  We only consider $s=0$. According to the method of   coupling by change of measures summarized in \cite[Section 1.1]{Wbook}, the main steps of the proof include:
\beg{enumerate}
\item[(S1)] Let $(X_t)_{t\ge 0}$ solve \eqref{E11} with $\L_{X_0}=\mu_0$. By the uniqueness we have   $\mu_t:=\P_t^*\mu_0= \L_{X_t}$, and  the equation \eqref{E11} reduces to
\beq\label{CP1} \d X_t= b_t(X_t,\mu_t)\d t+\si_t(X_t)\d W_t.\end{equation}
\item[(S2)] Construct a process $(Y_t)_{t\in [0,T]}$ such that for a weighted probability measure $\Q:=R_T\P,$
\beq\label{CP2} X_T=Y_T\ \Q\text{-a.s., \ \ and}\ \L_{Y_T}|_\Q=P_T^*\nu_0=:\nu_T. \end{equation}
\end{enumerate}
Obviously,  (S1) and (S2) imply
\beq\label{CP0} (P_Tf)(\mu_0)= \E[f(X_T)]\ \text{and}\  (P_T f)(\nu_0)= \E_\Q[f(Y_T)]=\E[R_Tf(X_T)],\ \ f\in \B_b(\R^d).\end{equation}
Combining this with   Young's inequality, we obtain the log-Harnack inequality:
\beq\label{LHdd}\beg{split} (P_T\log f)(\nu_0)&\le \E[R_T\log R_T]+ \log\E[f(X_T)]\\
&=\log (P_Tf)(\mu_0)+ \E[R_T\log R_T],\ \ f\in \B_b^+(\R^d).\end{split}\end{equation} \end{proof}

\subsubsection{The degenerate case}
  Consider the following distribution dependent stochastic Hamiltonian system for   $(X_t,Y_t)\in \R^{d_1}\times \R^{d_2}:$
\beq\label{E'} \beg{cases} \d X_t= \big(AX_t+BY_t)\d t,\\
\d Y_t= Z(t,(X_t,Y_t), \L_{(X_t,Y_t)}) \d t + \si_t \d W_t,\end{cases}\end{equation}
where $A$ is a $d_1\times d_1$-matrix, $B$ is a $d_1\times d_2$-matrix, $\si$ is a $d_2\times d_2$-matrix,
  $W_t$ is the $d_2$-dimensional Brownian motion on a complete filtration probability space $(\OO,\{\F_t\}_{t\ge 0}, \P)$, and
  $$Z: [0,\infty)\times \R^{d_1+d_2}\times \scr P_2(\R^{d_1+d_2})\to\R^{d_2},\ \ \si: [0,\infty)\to \R^{d_2}\otimes\R^{d_2}$$are measurable. We assume
\beg{enumerate} \item[$(H_4^2)$]  $\si(t)$ is invertible, there exists a locally bounded function $K: [0,\infty)\to [0,\infty)$ such that
$$\|\si(t)^{-1}\|\le K(t),\ \   |Z(t,x,\mu)-Z(t,y,\nu)|\le K(t)\big\{|x-y|+\W_2(\mu,\nu)\big\}$$ holds for all $ t\ge 0, \mu,\nu\in \scr P_2(\R^{d_1+d_2})$ and $x,y\in\R^{d_1+d_2}$, and $A,B$ satisfy the following Kalman's rank condition for some $k\ge 1$:
$${\rm Rank}[A^0B,\cdots, A^{k-1} B]=d_1, \ \ A^0:= I_{d_1\times d_1}. $$ \end{enumerate}
Obviously, this assumption implies $(H_3^1)$, so that \eqref{E'} has a unique solution $(X_t,Y_t)$ for any initial value $(X_0,Y_0)$ with
$\mu:=\L_{(X_0,Y_0)}\in \scr P_2(\R^{d_1+d_2}).$ Let $P^*_t\mu:=\L_{(X_t,Y_t)}$   and
$$(P_tf)(\mu) :=\int_{\R^{d_1+d_2} } f\d P_t^*\mu,\ \ t\ge 0, f\in \B_b(\R^{d_1+d_2}).$$ By \cite[Theorem 3.1]{W18}, the Lipschitz continuity of $Z$ implies
\beq\label{ES0} \W_2(P_t^*\mu, P_t^*\nu) \le \e^{Kt} \W_2(\mu,\nu),\ \ t\ge 0, \mu,\nu\in \scr P_2(\R^{d_1+d_2})\end{equation} for some constant $K>0.$
The following result is due to \cite[Section 5.1]{RW20b}.

\beg{thm}[\cite{RW20b}]\label{P1} Assume $(H_4^2)$. Then there exists an increasing  function $C: [0,\infty)\to (0,\infty)$ such that for any $T>0$,
\begin{equation*} (P_T\log f)(\nu)\le \log (P_T f)(\mu) + \ff {C(T)}{T^{4k-1}\land 1} \W_2(\mu,\nu)^2,\ \  \mu,\nu\in \scr P_2(\R^{d_1+d_2}), f\in \B^+_b(\R^{d_1+d_2}).\end{equation*}
Consequently,
\beq\label{LH'} \Ent(P_T^*\nu|P_T^*\mu) \le   \ff {C(T)}{T^{4k-1}\land 1}  \W_2(\mu,\nu)^2,\ \ T>0, \mu,\nu\in \scr P_2(\R^{d_1+d_2}).\end{equation}
 \end{thm}

\subsection{Bismut formula for the Lions derivative of $P_tf$}
We first introduce the intrinsic and Lions   derivatives for functionals of measures, then present the Bismut formula for the Lions derivative of $P_tf$ for non-degenerate and degenerate DDSDEs respectively.
The main results are taken from \cite{RW18}, see also \cite{BRW20} for  extensions to distribution-path dependent SDEs.

\subsubsection{Intrinsic and  Lions derivatives}

 \beg{defn} Let   $f: \scr P_2(\R^d)\to \R$.\beg{enumerate}
 \item[$(1)$]  If for any  $\phi\in L^2(\R^d\to\R^d;\mu)$,
  $$D_\phi^I f(\mu):= \lim_{\vv\downarrow 0} \ff{f(\mu\circ({\rm Id}+\vv \phi)^{-1})-f(\mu)}\vv\in\R$$ exists, and is a bounded linear functional in  $\phi$, we call   $f$ intrinsic differentiable  at   $\mu$. In this case, there exists a unique   $D^If(\mu)\in L^2(\R^d\to\R^d;\mu)$  such that
  $$\<D^If(\mu), \phi\>_{L^2(\mu)} = D^I_\phi f(\mu),\ \ \phi\in L^2(\R^d\to\R^d;\mu).$$
We call   $D^If(\mu)$ the  intrinsic derivative  of  $f$ at  $\mu$. If   $f$  is intrinsic differentiable at all $\mu\in \scr P_2(\R^d)$, we call it intrinsic differentiable on  $\scr P_2(\R^d)$ and denote
$$\|D^If(\mu)\|:= \|D^If(\mu)\|_{L^2(\mu)}=\bigg(\int_{\R^d} |D^If(\mu)|^2\d\mu\bigg)^{\ff 1 2}.$$
\item[$(2)$]     If    $f$ is intrinsic differentiable and  for any  $\mu\in\scr P_2(\R^d)$,
  $$\lim_{\|\phi\|_{L^2(\mu)}\to 0} \ff{f(\mu\circ({\rm Id}+\phi)^{-1})-f(\mu)- D^I_\phi f(\mu)}{\|\phi\|_{L^2(\mu)}}=0,$$
we call   $f$  $L$-differentiable on   $\scr P_2(\R^d)$.  In this case,   $D^If(\mu)$ is also denoted by   $D^Lf(\mu)$, and is called the   $L$-derivative of $f$ at $\mu$. \end{enumerate}
\end{defn}
Intrinsic derivative was first introduced in \cite{AKR} in the configuration space over a Riemannian manifold, while the $L$-derivative appeared in the Lecture notes \cite{Card} for the study of mean field games and is also called Lions derivative in references.

Note that    the derivative  $D^I f(\mu)\in L^2(\R^d\to\R^d;\mu)$  is   $\mu$-a.e.  defined. In applications, we take its continuous version  if exists.
 The following classes of $L$-differentiable functions are often used in analysis: \beg{enumerate}
 \item[{\bf (a)}]   $f\in C^{1}(\scr P_2(\R^d)):$    if  $f$ is   $L$-differentiable such that for every  $\mu\in\scr P_2(\R^d)$,
there exists a   $\mu$-version   $D^L f(\mu)(\cdot)$ such that    $D^L f(\mu)(x)$ is jointly continuous in    $(x,\mu)\in\R^d\times \scr P_2(\R^d)$.
 \item[{\bf (b)}]    $f\in C_b^1(\scr P_2(\R^d)):$   if   $f\in C^1(\scr P_2(\R^d))$ and
  $D^L f(\mu)(x)$   is bounded.
 \item[{\bf (c)}]   $f\in C^{2}(\scr P_2(\R^d)):$  if   $f\in C^{1}(\scr P_2(\R^d))$ and  $Df(\mu)(x)$ is  $L$-differentiable in  $\mu$ and differentiable in   $x\in\R^d$, such that  $\nn \{D^L f(\mu)\}(x)$  and   $$ (D^L)^2f(\mu)(x,y) := \big(\big\{D^L[D^Lf(\mu)(x)]_i (y)\big\}_j\big)_{1\le i,j\le d}\in \R^{d}\otimes\R^d$$  are jointly continuous in  $(\mu,x,y)\in \scr P_2(\R^d)\times\R^d\times\R^d.$
 \item[{\bf (d)}]   $f\in C_b^2(\scr P_2(\R^d)): $   if $f\in C^2(\scr P_2(\R^d))$ and
  all   derivatives  $D^Lf(\mu)(x), (D^L)^2f(\mu)(x,y)$ and $\nn (D^L f(\mu))(x)$ are bounded.
 \item[{\bf (e)}]   $f\in C^{1,1}(\R^d\times\scr P_2(\R^d)):$ if  $f$ is a continuous function on   $\R^d\times \scr P_2(\R^d)$  such that
  $f(\cdot,\mu)\in C^1(\R^d)$ for $\mu\in\scr P_2(\R^d)$,   $f(x,\cdot)\in C^{1}(\scr P_2(\R^d))$ for  $x\in \R^d$,
and
  $$ \nn f(x,\mu),\ \  \ D^Lf(x,\mu)(y)  $$ are  jointly
continuous in   $(x,\mu,y)\in \R^d\times \scr P_2(\R^d)\times \R^d.$
  If moreover these derivatives are bounded, we denote
  $f\in C_b^{1,1}(\R^d\times\scr P_2(\R^d)).$
 \item[{\bf (f)}]  {$f\in C^{2,2}(\R^d\times\scr P_2(\R^d))$}, if     $f$ is a continuous function on   $\R^d\times \scr P_2(\R^d)$  such that
 $f(\cdot,\mu)\in C^2(\R^d)$  for $\mu\in\scr P_2(\R^d)$,  $f(x,\cdot)\in C^{2}(\scr P_2(\R^d))$ for $x\in \R^d$,
  $$(D^L\nn f)(x,\mu)(y):= \big(\big\{D^L[\pp_{x_i} f(x,\mu)]\big\}_j\big)_{1\le i,j\le d}\in \R^d\otimes\R^d$$  exists, and all derivatives
  \beg{align*}& \nn f(x,\mu),\ \nn^2f(x,\mu),\ D^Lf(x,\mu)(y),\ (D^L\nn f)(x,\mu)(y)\\
&   \nn\{D^Lf(x,\mu)(\cdot)\}(y),\ (D^L)^2f(x,\mu)(y,z)\end{align*}  are  jointly
continuous in    $(x,\mu,y,z)\in \R^d\times \scr P_2(\R^d)\times \R^d\times \R^d.$
 If moreover these derivatives are bounded, we denote
 $f\in C_b^{2,2}(\R^d\times\scr P_2(\R^d)).$ \end{enumerate}

  Consider  $f\in \F C_b^2(\scr P_2(\R^d))$, i.e.      $$f(\mu)=g(\mu(h_1),\cdots, \mu(h_n)),\ n\ge 1, g\in C^2(\R^n), h_i\in C_b^2(\R^d).$$ Then it is easy to see that
 $f\in C_b^{2}(\scr P_2(\R^d))$  with
    \beg{align*}& D^L f(\mu)(y)=\sum_{i=1}^d (\pp_i g ) (\mu(h_1),\cdots, \mu(h_n))\nn h_i(y),\\
&  \nn \{D^L f(\mu)\}(y)=\sum_{i=1}^d (\pp_i g ) (\mu(h_1),\cdots, \mu(h_n))\nn^2 h_i(y),\\
& (D^L)^2f(\mu)(y,z)
 =\sum_{i,j=1}^d (\pp_i \pp_j g) (\mu(h_1),\cdots, \mu(h_n)) \{\nn h_i(y)\} \otimes \{\nn h_j (z)\}.\end{align*}

\subsubsection{ Bismut formula for non-degenerate DDSDEs}
Consider the DDSDE \eqref{E11} with coefficients satisfying the following assumption.
\beg{enumerate}\item[$(H_4^3)$] In addition to $(H_4^1)$,   $b_t,\si_t\in C^{1,1}(\R^d\times \scr P_2(\R^d))$ such  that
\beg{align*} &\max\Big\{\|\nn b_t(\cdot,\mu)(x)\|,  \|D^L b_t(x,\cdot)(\mu)\|,\  \ff 1 2\|\nn \si_t(\cdot,\mu)(x)\|^2,  \ff 1 2 \|D^L \si_t(x,\cdot)(\mu)\|^2  \Big\}\\
&\le K(t),\ \ \  t\ge 0, x\in\R^d, \mu\in \scr P_2(\R^d)   \end{align*}  holds  for some   continuous function $K:[0,\infty)\to [0,\infty).$     \end{enumerate}
By Theorem \ref{T1.1dd}, for any initial value $X_0\in L^2(\OO\to\R^d,\F_0,\P)$, \eqref{E11} has a unique solution $(X_t)_{t\ge 0}.$ Let $P_t^*\mu=\L_{X_t}$ for $\L_{X_0}=\mu$, and consider the $L$-derivative of the functionals in $\mu$:
  $$  P_Tf(\mu) := \E^\mu f(X_t)=\int_{\R^d} f(y) (P_T^*\mu )(\d y),\ \ T>0, f\in \B_b(\R^d).$$
Given $\phi\in L^2(\R^d\to\R^d,\mu)$,   the following  linear SDE has a unique solution   $v_t^\phi$ on $\R^d$:
\beq\label{DR0} \beg{split} \d v_t^\phi & =\Big\{\nn_{v_t^\phi} b(t,\cdot, \L_{X_t})(X_t)+ \big(\E \< D^L b(t,y,\cdot)(\L_{X_t})(X_t), v_t^\phi\>\big)\big|_{y=X_t}\Big\}\d t\\
&\quad + \big\{\nn_{v_t^\phi} \si(t,\cdot)(X_t)  \big\}\d W_t,\ \
v_0^\phi=\phi(X_0),\ t\ge 0.\end{split}\end{equation}
The following result is taken from \cite[Theorem 2.1 and Corollary 2.2]{RW19}.

\beg{thm}[\cite{RW19}] \label{TRW19-3.1} Assume  $(H_4^3)$.
Then for any   $f\in \B_b(\R^d),\mu\in\scr P_2(\R^d)$ and $T>0$, $P_Tf$ is  $L$-differentiable at $\mu$ such that  for any $g\in C^1([0,T])$ with $g_0=0$ and $g_T=1$,
$$    D_\phi^L(P_Tf)(\mu)  = \E\bigg[f(X_T)   \int_0^T \big\<  g_t'\si_t(X_t, \L_{X_t})^{-1} v_t^\phi,\ \d W_t\big\>\bigg], \  \phi\in L^2(\R^d\to\R^d,\mu), $$ where   $X_t$ solves $\eqref{E11}$ for $\L_{X_0}=\mu$. Moreover, the limit
\beq\label{*D*}  D_\phi^LP_T^*\mu:= \lim_{\vv\downarrow 0} \ff{P_T^*\mu\circ({\rm Id}+\vv\phi)^{-1}- P_T^*\mu}\vv= \psi P_T^*\mu\end{equation}   exists in the total variational norm, where
  $\psi$ is  the unique element in $L^2(\R^d\to\R, P_T^*\mu) $ such that
$\psi(X_T)=  \E\big(\int_0^T \big\<  g_t'\si_t(X_t, \L_{X_t})^{-1} v_t^\phi,\ \d W_t\big\>\big|X_T\big),$ and $(\psi P_T^*\mu)(A):= \int_A \psi\d P_T^*\mu,\ A\in\B(\R^d)$.
Consequently, for any $T>0, f\in \B_b(\R^d)$ and $\mu,\nu\in \scr P_2(\R^d),$
  \beg{align*} &\| D^L (P_Tf)(\mu)\|^2 \le \ff{(P_Tf^2)(\mu)-(P_t f(\mu))^2}{\int_0^T \lambda_t^{-2} \e^{-8K(t)t}\d t},\\
  &   \|P_T^*\mu-P_T^*\nu\|_{TV}^2:= 4\sup_{A\in \B(\R^d)} |(P_T^*\mu)(A)-(P_T^*\nu)(A)|^2
 \le \ff{4 \W_2(\mu,\nu)^2}{\int_0^T \lambda_t^{-2} \e^{-8K(t)t}\d t}.\end{align*}
  \end{thm}

\subsubsection{ Bismut formula for  degenerate DDSDEs}
Consider the following distribution dependent stochastic Hamiltonian system for $X_t=(X_t^{(1)}, X_t^{(2)})$ on $\R^{d_1+d_2}=\R^{d_1}\times \R^{d_2}$:
\beq\label{E5} \beg{cases} \d X_t^{(1)}= b^{(1)}_t(X_t)\d t,\\
\d X_t^{(2)} = b_t^{(2)}(X_t, \L_{X_t})\d t +\si_t \d W_t,\end{cases}\end{equation} where
  $(W_t)_{t\ge 0}$ is a $d_2$-dimensional Brownian motion as before, and  for each $t\ge 0$, $\si_t$ is an invertible $d_2\times d_2$-matrix,
$$b_t= (b_t^{(1)}, b_t^{(2)}): \R^{d_1+d_2}\times \scr P_2(\R^{d_1+d_2}) \to  \R^{d_1+d_2}$$ is measurable with $b^{(1)}_t(x,\mu)= b_t^{(1)}(x)$ independent of the distribution $\mu$.
Let $\nn=(\nn^{(1)}, \nn^{(2)})$ be the gradient operator on $\R^{d_1+d_2}=\R^{d_1}\times\R^{d_2}$, where $\nn^{(i)}$ is the gradient in the $i$-th component, $i=1,2$.
Let $\nn^2=\nn\nn$ denote the Hessian operator on $\R^{d_1+d_2}$. We assume

\beg{enumerate}\item[$(H_4^4)$] For every $t\ge 0$, $b_t^{(1)} \in C^2_b(\R^{d_1+d_2}\to\R^{d_1}),$ $b_t^{(2)}\in C^{1,1}(\R^{d_1+d_2}\times \scr P_2(\R^{d_1+d_2})\to\R^{d_2})$, and there exists an increasing function $K: [0,\infty)\to [0,\infty)$ such that
 $$\|\nn b_t(\cdot,\mu)(x)\|+ \|D^L b_t^{(2)}(x,\cdot)(\mu)\|  +\|\nn^2 b_t^{(1)}(x)\|  \le K(t),\ \
  t\ge 0, (x,\mu)\in\R^d\times \scr P_2(\R^d).$$
  There exist   $B\in \B_b([0,T]\to \R^{d_1}\otimes \R^{d_2}) $,  an increasing function $\theta\in C([0,T];\R^1)$ with $\theta_t>0$ for $t\in (0,T]$,  and $\vv\in (0,1)$  such that
  \beg{align*} & \<(\nn^{(2)} b^{(1)}_t - B_t)B_t^* a,a\>\ge -\vv |B_t^*a|^2,\ \ a\in \R^{d_1},\\
 & \int_0^t s(T-s) K_{T,s}B_sB_s^*K_{T,s}^*\d s\ge \theta_t  I_{d_1\times d_1},\ \ t\in (0,T],\end{align*}
 where for any $s\ge 0,$   $\{K_{t,s}\}_{t\ge s}$ is the unique solution of the following linear random ODE on $\R^{d_1}\otimes\R^{d_1}$:
$$
\ff{\d}{\d t}K_{t,s}=  (\nn^{(1)}b^{(1)}_t)(X_t) K_{t,s},\ \ \ t\ge s,  K_{s,s}=I_{d_1\times d_1}.$$
\end{enumerate}
\paragraph{Example 4.1.}  Let  $$b_t^{(1)}(x)= A x^{(1)} + Bx^{(2)},\ \  x=(x^{(1)}, x^{(2)}) \in \R^{d_1+d_2}$$ for some  $d_1\times d_1$-matrix $A$ and $d_1\times d_2$-matrix $B$. If
 the Kalman's rank condition
$$ \text{Rank}[B, AB, \cdots, A^{k}B]=d_1 $$ holds for some $k\ge 1$, then  $(H_4^4)$ is satisfied  with $\theta_t=c_T t$ for some constant $c_T>0$.

 \

According to the proof of \cite[Theorem 1.1]{WZ13}, $(H_4^4)$ implies that  the matrices
$$Q_t:=\int_0^t s(T-s) K_{T,s}\nn^{(2)}b_s^{(1)}(X_s) B_s^* K_{T,s}^* \d s,\ \   t\in (0,T]$$ are invertible with
\beq\label{Q}\|Q_t^{-1}\|\le\ff 1 { (1-\vv)\theta_t },\ \ t\in (0,T].\end{equation} For $(X_t)_{t\in [0,T]}$ solving \eqref{E5} with $\L_{X_0}=\mu\in \scr P_2(\R^{d_1+d_2})$ and $\phi=(\phi^{(1)},\phi^{(2)})\in L^2(\R^{d_1+d_2}\to\R^{d_1+d_2},\mu)$, let
\beg{align*} & \aa_t^{(2)}  =  \ff{T-t} T \phi^{(2)} (X_0) -\ff{t(T-t) B_t^*K_{T,t}^*}{\int_0^T \theta_s^2 \d s} \int_t^T  \theta_s^2 Q_s^{-1} K_{T,0}\phi^{(1)}(X_0)\d s\\
&\qquad\qquad -t(T-t) B_t^* K_{T,t}^*Q_T^{-1}\int_0^T\ff{T-s}TK_{T,s}\nn^{(2)}_{\phi^{(2)} (X_0)} b^{(1)}_s(X_s)\d s,\\
& \aa^{(1)}_t= K_{t,0}\phi^{(1)} (X_0)+\int_0^tK_{t,s}\nn^{(2)}_{\aa_s^{(2)}} b_s^{(1)}(X_s(x))\,\d s, \ \ t\in [0,T], \end{align*}  and define
\beq\label{B00}\beg{split}  h_t^\aa:= \int_0^t \si_s^{-1}  \Big\{ &\big(\E\<D^L b_s^{(2)}(y,\cdot)(\L_{X_s})(X_s), \aa_s\>\big)\big|_{y=X_s}\\
&\quad +\nn_{\aa_s} b_s^{(2)}(\cdot,\L_{X_s})(X_s) -(\aa_s^{(2)})'\Big\}\d s,\ \ t\in [0,T].\end{split}\end{equation}
Let $ (D^*, \D(D^*))$ be the   Malliavin divergence operator associated with the Brownian motion $(W_t)_{t\in [0,T]}$.  The following result is due to \cite[Theorem 2.3]{RW19}.

\beg{thm}[\cite{RW19}] \label{TRW19-4.2} Assume$(H_4^4)$. Then $h^\aa\in \D(D^*)$ with    $\E|D^*(h^{\aa})|^p<\infty$ for all $p\in [1,\infty)$. Moreover, for any $f\in \B_b(\R^{d_1+d_2})$ and $T>0$, $P_Tf$ is $L$-differentiable  such that
$$ D_\phi^L(P_Tf)(\mu)= \E\big[f(X_T)\, D^*(h^\aa)\big] $$ holds for $\mu\in \scr P_2(\R^{d_1+d_2}), \phi\in L^2(\R^{d_1+d_2}\to \R^{d_1+d_2},\mu)$  and $h^\aa$
in \eqref{B00}. Consequently:
\beg{enumerate} \item[$(1)$]  The formula $\eqref{*D*}$ holds for the unique $\psi\in L^2(\R^{d_1+d_2}\to\R, P_T^*\mu)$ such that $\psi(X_T)= \E(D^*(h^\aa)|X_T).$
\item[$(2)$] There exists  a  constant $c\ge 0$ such that for any $T>0$,
\beg{align*} &   \|D^L  (P_T f)(\mu)\|\le c\ss{ P_T|f|^2(\mu) -(P_Tf)^2(\mu)}  \ff{\ss{T}(T^2 +\theta_T )}{\int_0^{T} \theta_s^2\d s},\ \ f\in \B_b(\R^{d_1+d_2}),\\
& \|P_T^*\mu- P_T^*\nu\|_{TV} \le c \W_2(\mu,\nu) \ff{\ss{T}(T^2 +\theta_T) }{\int_0^{T} \theta_s^2\d s},\ \  \mu,\nu\in \scr P_2(\R^{d_1+d_2}).\end{align*}\end{enumerate}
 \end{thm}

\subsection{Lions derivative estimates on $P_tf$}
In this part we estimate $D^LP_Tf$ for DDSDE with    $\si$ also depending on $\mu$, which  thus extends  the corresponding derivative estimate presented in
Theorem \ref{TRW19-3.1}.

Consider the  DDSDE \eqref{DDSDE} with coefficients satisfying the following assumption which, by Theorem \ref{T1.1dd}, implies the well-posedness.
\beg{enumerate} \item[$(H_4^5)$] For any $t\ge 0$, $b_t,\si_t\in C^{1,1}(\R^d\times\scr P_2),$  and  there exists an increasing function $K: [0,\infty)\to [1,\infty)$ such that
 for any  $ t\ge 0, x,y\in\R^d$ and $ \mu\in \scr P_2(\R^d),$
$$  K_t^{-1} I_{d\times d} \le  (\si_t\si_t^*) (x,\mu)\le  K_t I_{d\times d},$$
\beg{align*}  &\ |b_t(x,\mu)|+ \|\nn b_t(\cdot,\mu)(x)\|+ \|D^L \{b_t(x,\cdot)\}(\mu)\| \\
&\qquad  +  \|\nn \{\si_t(\cdot,\mu)\}(x)\|^2+  \|D^L\{\si_t(x,\cdot)\}(\mu)\|^2
  \le K_t,\end{align*}
\beg{align*}   &\|D^L \{b_t(x,\cdot)\}(\mu)- D^L \{b_t(y,\cdot)\}(\mu)\|+ \|D^L \{\si_t(x,\cdot)\}(\mu)- D^L\{ \si_t(y,\cdot)\}(\mu)\|\\
&\le K_t |x-y|.
  \end{align*}
   \end{enumerate}
Let $P_{s,t}f(\mu):=\E[f(X_{s,t})]$ for $f\in \B_b(\R^d)$ and $(X_{s,t})_{t\ge s\geq 0}$ solving \eqref{DDSDE} with $\L_{X_{s,s}}=\mu\in \scr P_2(\R^d)$. The following result is due to \cite[Theorem 1.1]{HW20b}.

 \beg{thm}[\cite{HW20b}]  Assume $(H_4^5)$.  Then for any $t>s\geq 0$ and $f\in \B_b(\R^d)$, $P_{s,t}f$ is
  $L$-differentiable, and there exists an increasing function $C: [0,\infty)\to (0,\infty)$ such that
$$\|D^L P_{t} f(\mu)\|\le \ff{C_t \|f\|_\infty}{\ss {t-s}},\ \ t> s, f\in \B_b(\R^d).$$
 Consequently, for any $t> 0$ and $   \mu,\nu\in \scr P_2(\R^d),$
$$ \|P_{s,t}^*\mu- P_{s,t}^*\nu\|_{TV}:= 2\sup_{\|f\|_\infty\le 1} |P_{s,t}f(\mu)-P_{s,t}f(\nu)|\le \ff{2C_t}{\ss {t-s}} \W_2(\mu,\nu).$$
\end{thm}

\section{Exponential ergodicity in entropy}

 The convergence in entropy for stochastic systems is an important topic in both probability theory and  mathematical physics, and has been well studied for Markov processes by using the log-Sobolev inequality, see for instance \cite{BGL} and references therein. However, the existing results derived in the literature do not apply to DDSDEs.
   In 2003, Carrillo, McCann and Villani \cite{CMV} proved  the exponential convergence in a mean field entropy  of the  following granular media equation for   probability density functions $(\rr_t)_{t\ge 0}$ on $\R^d$:
 \beq\label{E0} \pp_t \rr_t= \DD\rr_t  + {\rm div} \big\{\rr_t\nn (V + W*\rr_t)\big\},\end{equation}
 where the internal potential $V\in C^2( \R^d)$  satisfies   $\Hess_V\ge \lambda  I_{d\times d}$ for a constant $\lambda>0$ and the $d\times d$-unit matrix $I_{d\times d}$,  and the interaction potential $W\in C^2(\R^d)$ satisfies
 $W(-x)=W(x)$ and  $\Hess_W\ge -\dd I_{d\times d}$
 for some constant $\dd\in [0, \lambda/2)$.   Recall  that  we write $M\ge \lambda I_d$ for a constant $\lambda$ and a $d\times d$-matrix $M$, if $\<Mv,v\>\ge \lambda |v|^2$ holds for any $v\in \R^d$.
 To introduce the mean field entropy, let  $\mu_V(\d x):= \ff{\e^{-V(x)}\d x}  {\int_{\R^d}\e^{-V(x)}\d x}$,   recall  the classical relative entropy
 $${\rm Ent}(\nu|\mu) := \beg{cases} \mu(\rr\log\rr),    &\text{if} \ \nu=\rr\mu,\\
  \infty, &\text{otherwise}\end{cases}$$ for $\mu,\nu\in \scr P(\R^d)$, and consider  the free energy functional
  $$E^{V,W}(\mu):= {\rm Ent}(\mu|\mu_V)+ \ff 1 2 \int_{\R^d\times\R^d} W(x-y) \mu(\d x)\mu(\d y),\ \ \mu\in \scr P(\R^d),$$
  where we set $E^{V,W}(\mu)=\infty$ if either ${\rm Ent}(\mu|\mu_V)=\infty$ or the integral term is not well defined.  Then the
 associated mean field entropy ${\rm Ent}^{V,W}$
   is defined by
\beq\label{ETP}  {\rm Ent}^{V,W}(\mu):= E^{V,W}(\mu) -  \inf_{\nu\in\scr P}  E^{V,W}(\nu),\ \ \mu\in \scr P(\R^d).\end{equation}
According to \cite{CMV},  for $V$ and $W$ satisfying the above mentioned conditions,   $E^{V,W}$ has a unique minimizer $\mu_\infty,$ and  $\mu_t(\d x):= \rr_t(x)\d x$  for probability density $\rr_t$ solving \eqref{E0} converges to $\mu_\infty$  exponentially    in the mean field  entropy:
$${\rm Ent}^{V,W}(\mu_t)\le \e^{-(\lambda-2\dd)t} {\rm Ent}^{V,W}(\mu_0),\ \ t\ge 0.$$
Recently, this result was generalized  in \cite{GLW} by establishing the uniform log-Sobolev inequality for the associated mean field particle systems, such that   ${\rm Ent}^{V,W}(\mu_t)$ decays exponentially   for a class of non-convex $V\in C^2(\R^d)$ and $W\in C^2(\R^d\times \R^d)$, where  $W(x,y)=W(y,x)$ and  $\mu_t(\d x):=\rr_t(x)\d x$ for $\rr_t$ solving the nonlinear PDE
\beq\label{E00} \pp_t \rr_t= \DD\rr_t  + {\rm div} \big\{\rr_t\nn (V + W\circledast\rr_t)\big\},\end{equation}
where
\beq\label{AOO}  W\circledast\rr_t:=\int_{\R^d}W(\cdot,y) \rr_t(y)\d y.\end{equation}     In this case,   $\Ent^{V,W}$ is defined in \eqref{ETP}  for  the free energy functional
$$E^{V,W}(\mu):=  {\rm Ent}(\mu|\mu_V)+ \ff 1 2 \int_{\R^d\times\R^d} W(x,y) \mu(\d x)\mu(\d y),\ \ \mu\in \scr P(\R^d).$$
To study \eqref{E00} using probability methods, we consider the
     following DDSDE with initial distribution $\mu_0$:
\beq\label{E0d'} \d X_t= \ss 2 \d B_t-\nn \big\{ V+ W\circledast\L_{X_t}\big\}(X_t)\d t,\end{equation}
where $B_t$ is the $d$-dimensional Brownian motion,  $\L_{X_t}$ is the distribution of $X_t$, and
\beq\label{A01} (W\circledast\mu)(x):=\int_{\R^d} W(x,y) \mu(\d y),\ \ x\in \R^d, \mu\in \scr P(\R^d)\end{equation}   provided the integral exists.
Let $\rr_t(x)=\ff{(\L_{X_t})(\d x)}{\d x},\ t\ge 0.$ By It\^o's formula and the integration by parts formula,   we have
\beg{align*}& \ff{\d}{\d t} \int_{\R^d} ( \rr_t f)(x)\d x =\ff{\d }{\d t} \E [f(X_t)] =\E\big[\big(\DD-\nn V- \nn \{W\circledast\rr_t\}\big)f(X_t)\big]\\
&= \int_{\R^d} \rr_t(x) \big\{\DD f -\<\nn V+ \nn \{W\circledast\rr_t\}, \nn f\>  \big\}(x)\d x\\
&= \int_{\R^d} f(x) \{\DD\rr_t + {\rm div}[\rr_t\nn V+ \rr_t \nn (W\circledast\rr_t)]\big\}(x)\d x,\ \ t\ge 0,\ f\in C_0^\infty(\R^d).\end{align*}
Therefore, $\rr_t$ solves \eqref{E00}. On the other hand, by this fact and the uniqueness of \eqref{E00} and \eqref{E0d'}, if $\rr_t$ solves \eqref{E00} with $\mu_0(\d x):=\rr_0(x)\d x$,
then $\rr_t(x)\d x=\L_{X_t}(\d x)$ for $X_t$ solving \eqref{E0d'} with $\L_{X_0}=\mu_0.$

To extend the study of \cite{CMV, GLW}, we   investigate the exponential convergence in entropy for  the following DDSDE on $\R^d$:
\beq\label{E100} \d X_t=  b(X_t,\L_{X_t})\d t+\si(X_t)\d W_t,\end{equation}
where $W_t$ is the $m$-dimensional Brownian  motion on a complete filtration probability space $(\OO, \{\F_t\}_{t\ge 0}, \P)$,
$$\si: \R^d\to \R^{d}\otimes \R^m,\ \ b:\R^d\times \scr P_2(\R^d)\to \R^d$$ are measurable.

Unlike in \cite{CMV, GLW} where the mean field particle systems are used to estimate the mean field entropy, we use
the log-Harnack inequality introduced in \cite{W10, RW10} and the Talagrand inequality developed in \cite{TAL, BGL, OV}.
Since the log-Harnack inequality is not yet available when $\si$  depends on the distribution,   in \eqref{E100} we only consider  distribution-free $\si$.

In the following subsections, we first present a criterion on the exponential convergence for DDSDEs by using the log-Harnack and Talagrand  inequalities,  then prove the exponential convergence for   granular media type equations  which generalizes  the framework of \cite{GLW}, and finally  consider exponential convergence for \eqref{E100} with non-degenerate and degenerate noises respectively.

\subsection{A criterion with  application to Granular media type equations}

In general, we consider the following DDSDE:
\beq\label{EP}  \d X_t= \si(X_t,\L_{X_t})\d W_t+ b(X_t,\L_{X_t})\d t,\end{equation}
where $W_t$ is the $m$-dimensional Brownian motion and
$$\si: \R^d\times \scr P_2(\R^d)\to \R^{d}\otimes \R^m,\ \ b:\R^d\times\scr P_2(\R^d) \to \R^d$$ are measurable.  We assume that this SDE is strongly and weakly well-posed
for square integrable initial values. It is in particular the case if $b$ is continuous on $\R^d\times\scr P_2(\R^d)$ and there exists a constant $K>0$ such that
\beq\label{KK00}\beg{split} &\<b(x,\mu)-b(y,\nu), x-y\>^+ +\|\si(x,\mu)-\si(y,\nu)\|^2\le K \big\{|x-y|^2 +\W_2(\mu,\nu)^2\},\\
& |b(0,\mu)|\le K\Big(1+\ss{\mu(|\cdot|^2)}\Big),\ \ x,y\in \R^d, \mu,\nu\in \scr P_2(\R^d),\end{split}\end{equation}
see for instance \cite{W18}. See also \cite{HW20c,Zhao} and references therein for the well-posedness of DDSDEs with singular coefficients. For any  $\mu\in \scr P_2(\R^d)$, let $P_t^*\mu=\L_{X_t}$ for the solution $X_t$ with initial distribution $\L_{X_0}=\mu$.  Let
$$P_t f(\mu)= \E[f(X_t)]=\int_{\R^d}f\d P_t^*\mu,\ \ t\ge 0, f\in \B_b(\R^d).$$
We have the following equivalence on the exponential convergence of $P_t^*\mu$ in $\Ent$ and $\W_2$.

\beg{thm}[\cite{RW20b}] \label{T0} Assume that $P_t^*$ has a unique invariant probability measure $\mu_\infty\in \scr P_2(\R^d)$ such that for some constants $t_0, c_0,C>0$ we have the log-Harnack inequality
\beq\label{LHI} P_{t_0} (\log f)(\nu)\le \log P_{t_0} f (\mu)+ c_0 \W_2(\mu,\nu)^2,\ \ \mu,\nu\in \scr P_2(\R^d),f\in\scr B^+_b(\R^d)\end{equation} and the Talagrand inequality
\beq\label{TLI} \W_2(\mu,\mu_\infty)^2\le C \Ent(\mu|\mu_\infty),\ \ \mu\in \scr P_2(\R^d).\end{equation}
\beg{enumerate} \item[$(1)$] If there exist constants $c_1,\lambda, t_1\ge 0$ such that
\beq\label{EW} \W_2(P_t^* \mu,\mu_\infty)^2
\le c_1\e^{-\lambda t} \W_2(\mu,\mu_\infty)^2,\ \ t\ge t_1, \mu\in \scr P_2(\R^d), \end{equation}
then
\beq\label{ET}\beg{split} & \max\big\{c_0^{-1} \Ent(P_t^*\mu|\mu_\infty), \W_2(P_t^*\mu, \mu_\infty)^2\big\} \\
&\le  c_1\e^{-\lambda (t-t_0)} \min\big\{\W_2(\mu,\mu_\infty)^2,  C\Ent(\mu|\mu_\infty)\big\},\ \ t\ge t_0+t_1, \mu\in \scr P_2(\R^d).\end{split}
 \end{equation}
\item[$(2)$]  If for some constants $\lambda, c_2,t_2>0$
\beq\label{EW'} \Ent (P_t^* \mu|\mu_\infty)\le c_2\e^{-\lambda t} \Ent(\mu|\mu_\infty),\ \ t\ge t_2, \mu\in \scr P_2(\R^d), \end{equation}
then
\beq\label{ET'}\beg{split} &\max\big\{\Ent(P_t^*\mu,\mu_\infty), C^{-1}  \W_2(P_t^* \mu,\mu_\infty)^2\big\}\\
&\le  c_2\e^{-\lambda (t-t_0)} \min\big\{c_0\W_2(\mu,\mu_\infty)^2, \Ent(\mu|\mu_\infty)\big\},\ \ t\ge t_0+t_2, \mu\in \scr P_2(\R^d).\end{split} \end{equation}
\end{enumerate}
  \end{thm}

When $\si\si^*$ is invertible and does not depend on the distribution, the log-Harnack inequality \eqref{LHI} has been established in \cite{W18}.
The Talagrand inequality was first found  in \cite{TAL} for $\mu_\infty$ being the Gaussian measure, and extended in \cite{BGL} to
$\mu_\infty$ satisfying the log-Sobolev inequality
\beq\label{LS0} \mu_\infty(f^2\log f^2)\le C \mu_\infty(|\nn f|^2),\ \ f\in C_b^1(\R^d), \mu_\infty(f^2)=1,\end{equation}
see  \cite{OV} for an earlier result under a curvature condition, and  see \cite{W04} for further extensions.

To illustrate this result, we consider   the  granular media type equation for probability density functions $(\rr_t)_{t\ge 0}$ on $\R^d$:
\beq\label{E01} \pp_t \rr_t=   {\rm div} \big\{ a\nn\rr_t + \rr_t a\nn (V + W\circledast\rr_t)\big\}, \end{equation}
 where   $W\circledast \rr_t$ is in \eqref{AOO}, and the functions
 $$ a: \R^d\to \R^d\otimes \R^d,\ \ V:\R^d\to\R,\ \  W:\R^d\times\R^d\to \R$$ satisfy the following assumptions.

\beg{enumerate} \item[$(H_5^1)$]  $a:=(a_{ij})_{1\le i,j\le d} \in C_b^2(\R^d\to \R^d\otimes\R^d)$,  and
$a\ge \lambda_aI_{d\times d}$ for some constant $\lambda_a>0$.
\item[$(H_5^2)$] $V\in C^2(\R^d), W\in C^2(\R^d\times\R^d)$ with  $W(x,y)=W(y,x)$, and there exist   constants $\kk_0\in\R$ and $\kk_1,\kk_2,\kk_0'>0$ such that
\beq\label{ICC1}  \Hess_V\ge \kk_0I_{d\times d},\ \ \kk_0' I_{2d\times 2d}\ge \Hess_{W}\ge \kk_0 I_{2d\times 2d},\end{equation}
\beq\label{ICC2} \<x, \nn V(x)\> \ge \kk_1|x|^2-\kk_2,\ \ x\in \R^d.\end{equation}
Moreover, for any $\lambda>0$,
\beq\label{ICC3} \int_{\R^d\times\R^d} \e^{-V(x)-V(y)-\lambda W(x,y)} \d x\d y<\infty.\end{equation}
\item[$(H_5^3)$] There exists a  function $b_0\in L^1_{loc}([0,\infty))$   with
$$r_0:=  \ff {\|\Hess_W\|_\infty} 4 \int_0^\infty \e^{\ff 1 4 \int_0^t b_0(s)\d s } \d t < 1$$
such that for any $x,y,z\in\R^d$,
\beg{align*}  \big\<y-x, \nn V(x)-\nn V(y) +\nn W(\cdot,z)(x)-\nn W(\cdot, z)(y)\big\>
  \le   |x-y| b_0(|x-y|).\end{align*} \end{enumerate}

\

 For any $N\ge 2$, consider the Hamiltonian for the   system of $N$ particles:
  $$H_N(x_1,\cdots, x_N)=\sum_{i=1}^N  V(x_i)+ \ff 1 {N-1} \sum_{1\le i<j\le N}^N W(x_i, x_j),$$
  and the corresponding finite-dimensional Gibbs measure
 $$\mu^{(N)}(\d x_1,\cdots, \d x_N)= \ff 1 {Z_N} \e^{-H_N(x_1,\cdots, x_N)}\d x_1\cdots\d x_N,$$
 where   $Z_N:=\int_{\R^{dN}} \e^{-H_N(x)}\d x <\infty$   due to \eqref{ICC3} in $(H_2)$.
 For any  $1\le i\le N$,  the  conditional marginal of   $\mu^{(N)}$  given $z\in \R^{d(N-1)}$ is  given by
\beg{align*} &\mu^{(N)}_z(\d x) := \ff 1 {Z_N(z)}  \e^{-H_N(x|z)} \d x,  \ \ Z_N(z):= \int_{\R^d} \e^{-H_N(x|z)} \d x,\\
&H_N(x|z):= V(x)-\log \int_{\R^{d(N-1)}} \e^{-\sum_{i=1}^{N-1} \{V(z_i) +\ff 1 {N-1} W(x, z_i)\}}\d z_1\cdots \d z_{N-1}.\end{align*}
We have the following result.

\beg{thm}[\cite{RW20b}] \label{C2.0} Assume $(H_5^1)$-$(H_5^3)$. If there is a constant $\beta>0$ such that the uniform   log-Sobolev inequality
\beq\label{LSS} \mu^{(N)}_z(f^2\log f^2)\le \ff 1 \bb \mu^{(N)}_z(|\nn f|^2),\ \ f\in C_b^1(\R^d), \mu^{(N)}_z(f^2)=1, N\ge 2, z\in \R^{d(N-1)}\end{equation} holds,
then there exists a unique   $\mu_\infty \in \scr P_2(\R^d)$ and a constant $c>0$ such that
\beq\label{WU} \W_2(\mu_t,\mu_\infty)^2 +\Ent(\mu_t|\mu_\infty)\le c \e^{-\lambda_a\bb (1-r_0)^2 t} \min\big\{\W_2(\mu_0, \mu_\infty)^2 +\Ent(\mu_0|\mu_\infty)\big\},\ \ t\ge 1\end{equation}
holds for   any probability density functions  $(\rr_t)_{t\ge 0}$ solving  \eqref{E01},  where $\mu_t(\d x) := \rr_t(x)\d x, t\ge 0.$  \end{thm}

This result    allows $V$ and W to be non-convex.   For instance, let $V=V_1+V_2\in C^2(\R^d)$ such that
$\| V_1\|_\infty\land\|\nn V_1\|_\infty<\infty$, $\Hess_{V_2}\ge \lambda I_{d\times d}$ for some $\lambda>0$, and   $W\in C^2(\R^d\times\R^d)$ with $\|W\|_\infty\land \|\nn W\|_\infty<\infty$. Then the uniform log-Sobolev inequality
\eqref{LSS} holds for some constant $\bb>0$.

\subsection{The non-degenerate case}
In this part, we make the following assumptions:
\beg{enumerate}  \item[$(H_5^4)$]  $b$ is continuous on $\R^d\times\scr P_2(\R^d)$ and there exists a constant $K>0$ such that $\eqref{KK00}$ holds.
\item[$(H_5^5)$]  $\si\si^*$ is invertible with $\lambda:=\|(\si\si^*)^{-1}\|_\infty<\infty$, and there exist constants $K_2>K_1\ge 0$ such that  for any $x,y\in \R^d$ and $\mu,\nu\in \scr P_2(\R^d)$,
$$ \|\si(x)-\si(y)\|_{HS}^2 + 2\<b(x,\mu)-b(y,\nu), x-y\>\le K_1 \W_2(\mu,\nu)^2 - K_2  |x-y|^2.$$
\end{enumerate}
According to  Theorem \ref{T1.1dd}, if $(H_3^1)$ holds, then for any initial value $X_0\in L^2(\OO\to\R^d,\F_0,\P)$, \eqref{E100}  has a unique solution which satisfies
$$\E\Big[\sup_{t\in [0,T]} |X_t|^2 \Big]<\infty,\ \ T\in (0,\infty).$$
Let $P_t^*\mu=\L_{X_t}$ for the solution $X_t$ with $\L_{X_0}=\mu.$
We have the following result.

\beg{thm}[\cite{RW20b}] \label{T1}  Assume  $(H_5^4)$ and $(H_5^5)$. Then  $P_t^*$ has a unique invariant probability measure $\mu_\infty$ such that
\beq\label{EX1} \max\big\{\W_2(P_t^*\mu,\mu_\infty)^2, \Ent(P_t^*\mu|\mu_\infty)\big\}\le  \ff{c_1}{t\land 1} \e^{-(K_2-K_1)t} \W_2(\mu,\mu_\infty)^2,\ \ t>0,  \mu\in \scr P_2(\R^d)\end{equation} holds for some constant $c_1>0$.
If moreover $\si\in C_b^2(\R^d\to \R^d\otimes \R^m)$, then there exists a constant $c_2>0$ such that  for any $\mu\in \scr P_2(\R^d), t\ge 1$,
\beq\label{EX2} \max\big\{\W_2(P_t^*\mu,\mu_\infty)^2, \Ent(P_t^*\mu|\mu_\infty)\big\} \le  c_2  \e^{-(K_2-K_1)t} \min\big\{\W_2(\mu,\mu_\infty)^2, \Ent(\mu|\mu_\infty)\big\}.    \end{equation}
\end{thm}

 \

 To illustrate this result, we consider the granular media equation \eqref{E00}, for which we take
\beq\label{SBB1} \si=\ss 2 I_{d\times d},~~~~~~ b(x,\mu)=  -\nn \big\{V+ W\circledast\mu\big\}(x),\ \ (x,\mu)\in \R^d\times \scr P_2(\R^d).\end{equation}
 The   following example is not included by   Theorem  \ref{C2.0} since the function $W$ may be non-symmetric.

\paragraph{Example 5.1 (Granular media equation).}  Consider \eqref{E00} with $V\in C^2(\R^d)$ and $W\in  C^2(\R^d\times\R^d)$ satisfying
\beq\label{CC1} \Hess_V\ge \lambda I_{d\times d},\ \ \Hess_W \ge \dd_1 I_{2d\times 2d},\ \ \|\Hess_W\|\le \dd_2\end{equation}
 for some constants $\lambda_1, \dd_2>0$ and $\dd_1\in \R$.
 If   $\lambda+\dd_1-\dd_2>0$, then there exists a unique    $\mu_\infty\in \scr P_2(\R^d)$ and a
 constant $c>0$ such that for any probability density functions   $(\rr_t)_{t\ge 0}$ solving \eqref{E00},   $\mu_t(\d x):=\rr_t(x)\d x$ satisfies
\beq\label{ERR0}\begin{split} &\max\big\{\W_2(\mu_t,\mu_\infty)^2,  \Ent(\mu_t|\mu_\infty)\big\} \\
&\le c\e^{-(\lambda+\dd_1-\dd_2)t} \min\big\{\W_2(\mu_0,\mu_\infty)^2,  \Ent(\mu_0|\mu_\infty)\big\},\ \ t\ge 1.
\end{split}\end{equation}

\beg{proof}  Let  $\si$ and $b$ be in \eqref{SBB1}.  Then \eqref{CC1} implies
 $(H_3^1)$   and
$$\<b(x,\mu)-b(y,\nu), x-y\> \le -(\lambda_1+\dd_1)|x-y|^2 + \dd_2 |x-y|\W_1(\mu,\nu),$$
where we have used the formula
$$\W_1(\mu,\nu)= \sup\{\mu(f)-\nu(f):\ \|\nn f\|_\infty\le 1\}.$$ So, by taking  $\aa=\ff{\dd_2}2$ and noting that $\W_1\le \W_2$, we obtain
\beg{align*} &\<b(x,\mu)-b(y,\nu), x-y\> \le -\big(\lambda+\dd_1-\aa \big)|x-y|^2 + \ff{\dd_2^2}{4\aa}  \W_1(\mu,\nu)^2\\
&\le - \Big(\lambda+\dd_1-\ff {\dd_2}2\Big) |x-y|^2+ \ff{\dd_2} 2  \W_2(\mu,\nu)^2,\ \ x,y\in\R^d,\mu,\nu\in \scr P_2(\R^d).\end{align*}
Therefore, if \eqref{CC1} holds for $\lambda+\dd_1-\dd_2>0$, Theorem \ref{T1} implies that $P_t^*$ has a unique invariant probability measure $\mu_\infty\in\scr P_2(\R^d)$, such that
\eqref{ERR0} holds for $\mu_0\in \scr P_2(\R^d)$.
  When $\mu_0\notin \scr P_2(\R^d)$, we have $\W_2(\mu_0,\mu_\infty)^2=\infty$ since $\mu_\infty\in \scr P_2(\R^d)$.
Combining this with the Talagrand inequality
$$\W_2(\mu_0,\mu_\infty)^2\le C \Ent(\mu_0|\mu_\infty)$$
for some constant $C>0$, see the proof of Theorem \ref{T1}, we have $\Ent(\mu_0|\mu_\infty)=\infty$ for $\mu_0\notin \scr P_2(\R^d)$, so that \eqref{ERR0} holds for all $\mu_0\in \scr P(\R^d)$.\end{proof}

 \subsection{The degenerate case}
When $\R^k$ with some $k\in \mathbb N$ is considered, to emphasize the space we use $\scr P(\R^k)$ ($\scr P_2(\R^k)$) to denote the class of probability measures (with finite second moment) on $\R^k$.
 Consider the following McKean-Vlasov stochastic Hamiltonian system for   $(X_t,Y_t)\in \R^{d_1+d_2}:= \R^{d_1}\times \R^{d_2}:$
\beq\label{E21} \beg{cases} \d X_t= BY_t\d t,\\
\d Y_t=   \ss 2 \d W_t - \Big\{ B^* \nn V(\cdot,\L_{(X_t,Y_t)})(X_t) + \bb B^* (BB^*)^{-1}X_t+Y_t\Big\}\d t,\end{cases}\end{equation}
where $\bb>0$ is a constant, $B$ is a $d_1\times d_2$-matrix such that $BB^*$ is invertible, and $$V: \R^{d_1}\times \scr P_2(\R^{d_1+d_2})\to \R^{d_2}$$ is measurable.
 Let
\beg{align*}& {\psi_B}((x,y),(\bar x,\bar y)):=\ss{|x-\bar x|^2 +|B(y-\bar y)|^2},\ \ (x,y), (\bar x, \bar y) \in \R^{d_1+d_2},\\
&\W_2^{\psi_B}(\mu,\nu):=\inf_{\pi\in \C(\mu,\nu)} \bigg\{ \int_{\R^{d_1+d_2}\times\R^{d_1+d_2}} {\psi_B}^2 \d\pi\bigg\}^{\ff 1 2},\ \ \mu,\nu\in\scr P_2(\R^{d_1+d_2}).\end{align*}
We assume
\beg{enumerate} \item[$(H_5^6)$]   $V(x,\mu)$ is differentiable in $x$ such that $\nn V(\cdot,\mu)(x)$ is Lipschitz continuous in $(x,\mu)\in \R^{d_1}\times \scr P_2(\R^{d_1+d_2}).$
Moreover, there exist   constants $\theta_1, \theta_2\in \R$ with
\beq\label{AH1} \theta_1+\theta_2<\bb,\end{equation}   such that  for any $(x,y), (x',y')\in \R^{d_1+d_2}$ and $\mu,\mu'\in \scr P_2(\R^{d_1+d_2})$,
\beq\label{AH2}\beg{split} & \big\< BB^* \{\nn V(\cdot,\mu)(x)-\nn V(\cdot,\mu')(x')\},  x-x'+(1+\bb)B(y-y')\big\>\\
&\ge -\theta_1{\psi_B} ((x,y), (x',y'))^2 -\theta_2 \W_2^{\psi_B} (\mu,\mu')^2.\end{split}\end{equation}
\end{enumerate}
Obviously, $(H_5^6)$ implies $(H_3^1)$ for $d=m=d_1+d_2$,  $\si= {\rm diag}\{0,\ss 2 I_{d_2\times d_2}\}$, and
$$b((x,y),\mu)= \big(By, - B^* \nn V(\cdot,\mu)(x) - \bb B^* (BB^*)^{-1}x-y\big).$$ So, according to \cite{W18},  \eqref{E21} is well-posed for  any initial value in $L^2(\OO\to \R^{d_1+d_2}, \F_0,\P)$.
Let $P_t^*\mu=\L_{(X_t,Y_t)}$ for the solution with initial distribution $\mu\in \scr P_2(\R^{d_1+d_2}).$

\beg{thm}[\cite{RW20b}] \label{T2} Assume $(H_5^6)$. Then $P_t^*$ has a unique invariant probability measure $\mu_\infty$ such that for any  $t>0$ and $\mu\in \scr P_2(\R^{d_1+d_2}),$
\beq\label{KK0} \max\big\{\W_2(P_t^*\mu,\mu_\infty)^2, \Ent(P_t^*\mu|\mu_\infty) \big\}\le \ff{c\e^{-2\kk t}}{(1\land t)^{3}}
\min\big\{\Ent(\mu|\mu_\infty), \W_2(\mu,\mu_\infty)^2\big\}\end{equation}
holds for some constant $c>0$ and
\beq\label{KK} \kk:=\ff{ 2(\bb-\theta_1-\theta_2)}{2+2\bb+\bb^2+\ss{\bb^4+4}}>0.\end{equation}
\end{thm}

\paragraph{Example 5.2 (Degenerate granular media equation).}  Let $ m\in \mathbb N$ and $W\in C^2(\R^m\times\R^{2m}).$ Consider  the following PDE for probability density functions $(\rr_t)_{t\ge 0}$ on $\R^{2m}$:
\beq\label{*EN2} \pp_t\rr_t(x,y)= \DD_y \rr_t(x,y) -\<\nn_x\rr_t(x,y), y\>+ \<\nn_y \rr_t(x,y), \nn_x (W\circledast\rr_t)(x) + \bb x+y\>,\end{equation}
where $\bb>0$ is a constant,   $\DD_y, \nn_x,\nn_y$ stand for the Laplacian in $y$ and the   gradient operators in $x,y$ respectively,  and
$$(W\circledast \rr_t)(x):=\int_{\R^{2m}} W(x,z)\rr_t(z) \d z,\ \ x\in \R^m.$$   If there exists a constant $\theta\in \big(0, \ff{2\bb}{1+3\ss{2+2\bb+\bb^2}}\big)$ such that
\beq\label{*EN3} |\nn W(\cdot,z)(x)-\nn W(\cdot,\bar z)(\bar x) |\le \theta \big(|x-\bar x|+|z-\bar z|\big),\ \ x,\bar x\in \R^m, z,\bar z\in \R^{2m},\end{equation}
then  there exists a unique probability measure  $\mu_\infty\in\scr P_2( \R^{2m})$   and a constant $c>0$ such that for any probability density functions $(\rr_t)_{t\ge 0}$ solving \eqref{*EN2},  $\mu_t(\d x):=\rr_t(x)\d x$ satisfies
\beq\label{ERR} \max\big\{\W_2(\mu_t,\mu_\infty)^2, \Ent(\mu_t|\mu_\infty) \big\}   \le c\e^{-\kk t} \min\big\{\W_2(\mu_0,\mu_\infty)^2, \Ent(\mu_0|\mu_\infty)\big\},\ \ t\ge 1\end{equation} holds for
$\kk= \ff{2\bb- \theta\big(1+ 3\ss{2+2\bb+\bb^2}\big)}{2+2\bb+\bb^2+\ss{\bb^4+4}}>0.$

\beg{proof} Let $d_1=d_2=m$ and $(X_t,Y_t)$ solve \eqref{E21} for
\beq\label{GPP} B:= I_{m\times m},\ \ V(x, \mu):=  \int_{\R^{2m}} W(x,z)\mu(\d z).\end{equation}
Let $\rr_t(z)= \ff{\L_{(X_t, Y_t)}(\d z)} {\d z}.$
By It\^o's formula and integration by parts formula, for any $f\in C_0^2(\R^{2m})$ we have
\beg{align*}& \ff{\d}{\d t} \int_{\R^{2m}} ( \rr_t f)(z)\d z =\ff{\d }{\d t} \E [f(X_t,Y_t)]\\
&= \int_{\R^{2m}}\rr_t(x,y)\big\{\DD_y f(x,y)+\<\nn_x f(x,y), y\> -\<\nn_y f(x,y), \nn_x V(x, \rr_t(z)\d z)+\bb x+y\>\big\}\d x\d y\\
&= \int_{\R^{2m}} f(x,y) \big\{\DD_y\rr_t (x,y)-  \<\nn_x \rr_t(x,y), y\> +\<\nn_y \rr_t(x,y), \nn_x \mu_t(W(x,\cdot)) +\bb x+y\>\big\}\d x\d y.
\end{align*} Then $\rr_t$ solves \eqref{*EN2}. On the other hand, by the uniqueness of of \eqref{E21} and \eqref{*EN2},  for any solution $\rr_t$ to \eqref{*EN2} with $\mu_0(\d z):= \rr_0(z)\d z\in \scr P_2( \R^{2m})$ for $d=2m$, $\rr_t(z)\d z=\L_{(X_t,Y_t)}(\d z)$ for the solution to \eqref{E21} with initial distribution $\mu_0$.
So, as explained in the proof of Example 2.1, by Theorem \ref{T2} we only need to verify $(H_5^6)$ for $B, V$ in \eqref{GPP}  and
\beq\label{TTH} \theta_1= \theta\Big(\ff 1 2 +\ss{2+2\bb+\bb^2}\Big),\ \ \theta_2= \ff \theta 2 \ss{2+2\bb+\bb^2},\end{equation}
so that the desired assertion holds for
$$ \kk:=\ff{ 2(\bb-\theta_1-\theta_2)}{2+2\bb+\bb^2+\ss{\bb^4+4}}=  \ff{2\bb- \theta(1+ 3\ss{2+2\bb+\bb^2})}{2+2\bb+\bb^2+\ss{\bb^4+4}}.$$
By \eqref{*EN3} and $V(x,\mu):=\mu(W(x,\cdot))$, for any constants $\aa_1,\aa_2,\aa_3>0$ we have
\beg{align*} I &:= \big\<\nn V(\cdot,\mu)(x)- \nn V(\cdot,\bar\mu)(\bar x), x-\bar x+(1+\bb)(y-\bar y)\big\> \\
&= \int_{\R^{2m}} \big\< \nn W(\cdot, z)(x)- \nn W(\cdot, z)(\bar x), x-\bar x +(1+\bb)(y-\bar y)\big\>\mu(\d z) \\
&\qquad +\big\<\mu(\nn_{\bar x}W(\bar x,\cdot))-\bar\mu(\nn_{\bar x} W(\bar x,\cdot)),  x-\bar x +(1+\bb)(y-\bar y)\big\>\\
&\ge -\theta\big\{ |x-\bar x|+ \W_1(\mu,\bar \mu)\big\} \cdot\big(|x-\bar x|+(1+\bb)|y-\bar y|\big)\\
&\ge - \theta (\aa_2+\aa_3) \W_2(\mu,\bar\mu)^2- \theta\Big\{\Big(1+\aa_1+\ff 1 {4\aa_2}\Big)|x-\bar x|^2 +(1+\bb)^2 \Big(\ff 1 {4 \aa_1}+\ff 1 {4\aa_3}\Big)|y-\bar y|^2\Big\}.\end{align*}
Take
$$\aa_1= \ff{\ss{2+2\bb+\bb^2}-1} 2,\ \ \aa_2= \ff 1 {2\ss{2+2\bb+\bb^2 }},\ \ \aa_3=\ff {(1+\bb)^2}{2\ss{2+2\bb+\bb^2 }}.$$
We have
\beg{align*}&1+\aa_1+\ff 1 {4\aa_2} =\ff 1 2+ \ss{2+2\bb+\bb^2},\\
&(1+\bb)^2\Big(\ff 1 {4 \aa_1}+\ff 1 {4\aa_3}\Big)=\ff 1 2+ \ss{2+2\bb+\bb^2},\\
& \ \aa_2+\aa_3= \ff 1 2 \ss{2+2\bb+\bb^2}.\end{align*}
Therefore,
 $$I\ge -\ff \theta 2 \ss{2+2\bb+\bb^2} \W_2(\mu,\bar\mu)^2 -\theta\Big(\ff 1 2 + \ss{2+2\bb+\bb^2}\Big) |(x,y)-(\bar x,\bar y)|^2,$$
 i.e.  $(H_5^6)$  holds for $B$ and $V$ in \eqref{GPP} where $B=I_{m\times m}$ implies that ${\psi_B}$ is the Euclidean distance on $\R^{2m}$, and for    $\theta_1,\theta_2$ in \eqref{TTH}.
\end{proof}

\section{Donsker-Varadhan large deviations}

The LDP (large deviation principle)  is a fundamental tool characterizing the asymptotic behaviour of probability measures  $\{\mu_\vv\}_{\vv>0}$ on a topological space $E$, see  \cite{DZ} and references within.
 Recall that  $\mu_\vv$ for small $\vv>0$ is said to
satisfy the LDP  with speed $\lambda(\vv)\to +\infty$ (as $\vv\to 0$) and rate
function $I: E\to [0,+\infty]$, if $I$ has compact level sets (i.e. $\{I\le r\}$ is compact  for $r\in {\mathbb R}^+$),
and for any Borel subset $A$ of $E$,
$$
 -\inf_{ A^o} I\le \liminf_{\vv\to 0}
   \frac 1{\lambda(\vv)} \log \mu_{\vv}(A)\le \limsup_{\vv\to 0}
   \frac 1{\lambda(\vv)} \log \mu_{\vv}(A)\le
  -\inf_{\bar A} I,
$$
where $A^o$ and $\bar A$ stand for    the
interior and the closure of $A$ in $E$ respectively.

In this part, we consider the  Donsker-Varadhan type long time LDP  \cite{DV}  for  $\mu_\vv:=\L_{L_{\vv^{-1}}},$  where
$$L_{t}:=\ff 1 t \int_0^{t} \dd_{X(s)}\d s,\ \ t>0$$  is the   empirical measure for a path-distribution dependent SPDE.

Let  $(\H, \<\cdot,\cdot\>, |\cdot|)$ be a separable Hilbert space. For a fixed constant $r_0>0$,
a  path $\xi\in \C:=C([-r_0,0];\H)$ stands for a sample of the history with  time length $r_0$. Recall that  $\C$ is a  Banach space   with the uniform norm
$$\|\xi\|_\infty:= \sup_{\theta \in [-r_0,0]}|\xi (\theta)|,\ \  \xi\in \C.$$
  For any map $\xi(\cdot): [-r_0,\infty)\to \H$ and any time $t\ge 0$, its segment $\xi_\cdot: [0,\infty)\to \C$ is defined by
      $$\xi_t(\theta):= \xi(t+\theta),\ \ \theta\in [-r_0,0], t\ge 0.$$
Let $\scr P(\C)$ denote the space of all probability measures on $\C$ equipped with the weak topology, and let
 $\L_{\eta}$ stand for the distribution of a random variable $\eta$.
Consider the following path-distribution dependent  SPDE on $\H$:
\begin{equation}\label{*1}
\d X(t)=\{AX(t)+b(X_t,\L_{X_t})\}\d t+\si(\L_{X_t})\d W(t),~~t\ge 0,
\end{equation}
where
\beg{enumerate} \item[$\bullet$]   $(A,\D(A))$ is a negative definite self-adjoint operator on $\H$;
\item[$\bullet$] $W(t)$ is the cylindrical Brownian motion on a separable Hilbert space $\tt\H$; i.e.
$$W(t)=\sum_{i=1}^\infty B_i(t)\tt e_i,\ \ t\ge 0$$ for an orthonormal basis $\{\tt e_i\}_{i\ge 1}$ on $\tt\H$ and a sequence of independent one-dimensional
 Brownian motions $\{B_i\}_{i\ge 1} $
on a complete filtration probability space $(\OO,\F,\{\F_t\}_{t\ge 0}, \P)$,
where $\F_0$ is rich enough such that for any $\pi\in\scr P(\C\times\C)$ there exists a  $\C\times\C$-valued random variable $\xi$ on $(\OO,\F_0,\P)$ such that $\L_{\xi}=\pi$.
\item[$\bullet$]  $b:\C\times\scr P(\C)\to\H,\ \ \si:  \scr P(\C)\to\mathbb L(\tt\H;\H)$ are
measurable. \end{enumerate}

 Let $X_t^\nu$ denote the mild segment solution with initial distribution  $\nu\in \scr P(\C)$, which is a continuous adapted process on $\C$.  We study the long time LDP for the empirical measure
$$L_t^\nu:= \ff 1 t \int_0^t \dd_{X_s^\nu}\d s,\ \ t>0.$$

 \beg{defn} \label{D2.1} Let $\scr P(\C)$ be equipped with the weak topology, let $\scr A\subset \scr P(\C)$, and let $J: \scr P(\C)\to [0,\infty]$ have compact level sets, i.e. $\{J\le r\}$ is compact in $\scr P(\C)$ for any $r>0$.
 \beg{enumerate} \item[(1)] $\{L_t^\nu\}_{\nu\in \scr A}$ is said to satisfy the upper bound uniform LDP with rate function $J$, denoted by $\{L_t^\nu\}_{\nu\in \scr A}\in LDP_u(J),$ if for any closed   $A\subset \scr P(\C),$
 $$\limsup_{t\to\infty} \ff 1 t\sup_{\nu\in \scr A} \log \P(L_t^\nu\in A)\le - \inf_AJ.$$
 \item[(2)] $\{L_t^\nu\}_{\nu\in \scr A}$ is said to satisfy the lower  bound uniform LDP with rate function $J$, denoted by $\{L_t^\nu\}_{\nu\in \scr A}\in LDP_l(J),$ if for any open  $A\subset \scr P(\C),$
 $$\liminf_{t\to\infty} \ff 1 t\inf_{\nu\in \scr A} \log \P(L_t^\nu\in A)\ge - \inf_AJ.$$
 \item[(3)] $\{L_t^\nu\}_{\nu\in \scr A}$ is said to satisfy the uniform LDP with rate function $J$, denoted by
 $\{L_t^\nu\}_{\nu\in \scr A}\in LDP(J),$ if $\{L_t^\nu\}_{\nu\in \scr A}\in LDP_u(J)$ and $\{L_t^\nu\}_{\nu\in \scr A}\in    LDP_l(J).$
\end{enumerate} \end{defn}

We investigate the long time LDP for \eqref{*1} in the following   three situations respectively:
 \beg{enumerate} \item[1)] $r_0=0$ and $\H$ is finite-dimensional;
 \item[2)]   $r_0=0$ and $\H$ is infinite-dimensional;   \item[3)] $r_0>0$ and $\si$ is constant.\end{enumerate}
 When $r_0>0$ and $\si$ is non-constant, the Donsker-Varadhan LDP is still unknown.

 To state establish the LDP, we recall   the Feller property, the strong Feller property and the irreducibility for a (sub-) Markov operator $P$.
 Let $\B_b(\C)$ (resp. $C_b(\C)$) be the space of bounded measurable (resp. continuous) real functions on $\C$.
Let $P$ be  a sub-Markov operator on $\B_b(\C)$, i.e. it is a positivity-preserving linear operator with $P1\le 1$.     $P$ is called strong Feller if $P\B_b(\C)\subset C_b(\C)$, is called Feller if $PC_b(\C)\subset C_b(\C)$, and is called $\mu$-irreducible for some $\mu\in \scr P(\C)$ if
$\mu(1_AP1_B)>0$ holds for any $A,B\in \B(\C)$ with $\mu(A)\mu(B)>0.$

  \subsection{ Distribution dependent SDE on $\R^d$ }

Let $r_0=0,$  $\H=\R^d$ and $\tt\H=\R^m$ for some $d,m\in \mathbb N$.
  In this case, we combine the linear term $Ax$ with the drift term $b(x,\mu)$,  so that \eqref{*1} reduces to
  \beq\label{E0dd} \d X(t)= b(X(t),\L_{X(t)})\d t + \si(\L_{X(t)}) \d W(t),\end{equation}
  where $b: \R^d\times \scr P_2(\R^d)\to \R^d,$   $\si:  \scr P_2(\R^d)\to \R^d\otimes\R^m $ and $W(t)$ is the $m$-dimensional Brownian motion. We assume
  \beg{enumerate}\item[$(H_6^1)$] $b$ is continuous, $\si$ is  bounded and   continuous such that
  $$2\<b(x,\mu) -b(y,\nu),x-y\>+ \|\si(\mu)-\si(\nu)\|_{HS}^2\le -\kk_1|x-y|^2 +\kk_2\W_2(\mu,\nu)^2$$ holds for some constants $\kk_1>\kk_2\ge 0$ and all $x,y\in\R^d, \mu,\nu\in \scr P_2(\R^d)$.
\end{enumerate}
Under $(H_6^1)$, for any $X(0)\in L^2(\OO\to\R^d,\F_0,\P)$,  the equation \eqref{E0dd} has a unique solution. We write $P_t^*\mu=\L_{X(t)}$ if $\L_{X(0)}=\mu$.  By \cite[Theorem 3.1(2)]{W18}, $P_t^*$ has a unique invariant probability measure $\bar\mu\in \scr P_2(\R^d)$ such that
\beq\label{EXPO}\W_2(P_t^*\nu, \bar\mu)^2\le\e^{-(\kk_1-\kk_2)t} \W_2(\nu,\bar\mu)^2,\ \ t\ge 0, \nu\in\scr P_2(\R^d).\end{equation}
Consider the reference SDE
\beq\label{E0'} \d \bar X(t)= b(\bar X(t),\bar\mu)\d t + \si(\bar\mu) \d W(t).\end{equation}
It is standard that under $(H_6^1)$ the equation \eqref{E0'} has a unique solution $\bar X^x(t)$ for any starting point $x\in \R^d,$ and $\bar\mu$ is the unique invariant probability measure of the associated Markov semigroup
$$\bar P_t f(x):= \E[f(\bar X^x(t))],\ \ t\ge 0, x\in \R^d, f\in \B_b(\R^d).$$
Consequently, $\bar P_t$ uniquely extends to $L^\infty(\bar\mu)$. If $f\in L^\infty(\bar\mu)$ satisfies
$$\bar P_t f= f+\int_0^t \bar P_s g\d s,\ \ \bar\mu\text{-a.e.}$$ for some $g\in L^\infty(\bar\mu)$ and  all $t\ge 0$, we write
$f\in \D(\bar{\scr A})$ and denote $\bar{\scr A} f=g$. Obviously, we have $\D(\bar{\scr A})\supset C_c^\infty(\R^d):=\{f\in C_b^\infty(\R^d): \nn f\text{\ has\ compact\ support}\}$ and
$$\bar{\scr A} f(x)= \ff 1 2 \sum_{i,j=1}^d \{\si\si^*\}_{ij}(\bar\mu)\pp_{i}\pp_j f(x)+ \sum_{i=1}^d b_i(x,\bar\mu) \pp_if(x),\ \ f\in C_c^\infty(\R^d).$$

The Donsker-Varadhan level 2 entropy function $J$ for the diffusion process generated by $\bar{\scr A}$ has compact level sets in $\scr P(\R^d)$ under the $\tau$ and weak topologies, and by \cite[3.11]{RW20a}, we have
$$ J(\nu)= \beg{cases} \sup\big\{\int_{\R^d} \ff{-\bar{\scr A} f}{f} \d\nu:\   1\le f\in \D(\bar{\scr A})\big\}, &\text{if}\ \nu\ll\mu,\\
 \infty, &\text{otherwise}.\end{cases} $$

\beg{thm}[\cite{RW20a}]  \label{T01} Assume $(H_6^1)$. For any $r,R>0$, let
$\scr B_{r,R}=\big\{\nu\in \scr P(\R^d): \nu(\e^{|\cdot|^r})\le R\big\}.$
\beg{enumerate}
\item[$(1)$] We have $\{L^\nu_t\}_{\nu\in \scr B_{r,R}}\in LDP_u(J)$  for all $r,R>0$.
If $\bar P_t$ is strong Feller and $\bar\mu$-irreducible for some $t>0$, then $\{L^\nu_t\}_{\nu\in \scr B_{r,R}}\in LDP(J)$  for all $r,R>0$.
\item[$(2)$] If there exist constants $\vv,c_1,c_2>0$ such that
 \beq\label{ABC}  \<x, b(x,\nu)\>  \le c_1-c_2|x|^{2+\vv},\ \ x\in \R^d,\nu\in \scr P_2(\R^d),\end{equation} then $\{L^\nu_t\}_{\nu\in \scr P_2(\R^d)} \in LDP_u(J)$. If moreover $\bar P_t$ is strong Feller and $\bar\mu$-irreducible  for some $t>0$, then $\{L^\nu_t\}_{\nu\in \scr P_2(\R^d)}\in LDP(J).$   \end{enumerate}
\end{thm}

To apply this result, we first recall some facts on the strong Feller property and the irreducibility of diffusion semigroups.
\begin{rem}\label{Remark 2.1.} (1) Let $\bar P_t$ be the (sub-)Markov semigroup generated by the second order differential operator
$$\bar{\scr A}:= \sum_{i=1}^m U_i^2 +U_0,$$
where $\{U_i\}_{i=1}^m$ are $C^1$-vector fields and $U_0$ is a continuous vector field. According to   \cite[Theorem 5.1]{LAP}, if       $\{U_i: 1\le i\le m\}$    together with their Lie brackets with $U_0$  span $\R^d$ at any point (i.e. the H\"ormander condition holds),  then  the Harnack inequality
$$ P_t f(x)\le \psi(t,s,x,y)   P_{t+s} f(y),\ \ t,s>0, x,y\in \R^d, f\in \B^+(\R^d)$$
for some map $\psi: (0,\infty)^2\times (\R^d)^2\to (0,\infty).$ Consequently, if moreover $\bar P_t$ has an invariant probability measure $\bar\mu$, then  $\bar P_t$ is $\bar\mu$-irreducible for any $t>0.$   Finally, if $\{U_i\}_{0\le i\le m}$ are smooth with bounded derivatives of all orders, then
the above H\"ormander condition implies that $\bar P_t$ has smooth heat kernel  with respect to the Lebesgue measure, in particular  it is strong Feller for any $t>0.$

(2) Let $\bar P_t$ be the  Markov semigroup generated by
$$\bar{\scr A}:= \sum_{i,j=1}^d \bar a_{ij} \pp_i\pp_j  +\sum_{i=1}^d \bar b_i \pp_j,$$
where $(\bar a_{ij}(x))$ is strictly positive definite for any $x$, $\bar a_{ij}\in H_{loc}^{p,1}(\d x)$ and $\bar b_i\in L_{loc}^p(\d x)$ for some $p>d$ and
all $1\le i,j\le d.$ Moreover, let $\bar\mu$ be an invariant probability measure of $\bar P_t$. Then by \cite[Theorem 4.1]{BKR},
$\bar P_t$ is strong Feller for all $t>0$. Moreover, as indicated in (1) that   \cite[Theorem 5.1]{LAP} ensures the $\bar\mu$-irreducibility of $\bar P_t$ for $t>0$.
\end{rem}
  \

We present below two examples to illustrate this result, where the first is a distribution dependent perturbation of the Ornstein-Ulenbeck process, and the second is the distribution dependent stochastic Hamiltonian system.

\paragraph{Example 6.1.} Let $\si(\nu)=I+ \vv\si_0(\nu)$ and $b(x,\nu)= -\ff 1 2 (\si\si^*)(\nu)x$, where $I$ is the identity matrix, $\vv>0$ and $\si_0$ is a bounded Lipschitz continuous  map from $\scr P_2(\R^d)$ to $\R^d\otimes\R^d$. When $\vv>0$ is small enough, assumption  {\bf $(H_1)$} holds and
that $\bar P_t$ satisfies conditions in Remark \ref{Remark 2.1.}(2).  So, Theorem \ref{T01}(1)   implies    $\{L^\nu_t\}_{\nu\in \scr B_{r,R}}\in LDP(J)$  for all $r,R>0$.

If we take $b(x,\nu)= - x- c |x|^\theta x$ for some constants
$c,\theta>0$, then  when $\vv>0$ is small enough such that  $(H_1)$ and \eqref{ABC} are satisfied, Theorem \ref{T01}(2) and  Remark \ref{Remark 2.1.}(2) imply $\{L^\nu_t\}_{\nu\in \scr P_2(\R^d)}\in LDP(J).$

\paragraph{Example 6.2.} Let $d=2m$ and consider the following distribution dependent SDE for $X(t)=(X^{(1)}(t), X^{(2)}(t))$ on $\R^{m}\times \R^m:$
$$\beg{cases} \d X^{(1)}(t)= \{X^{(2)}(t)- \lambda X^{(1)}(t)\}\d t\\
\d X^{(2)}(t)= \{Z(X(t),\L_{X(t)}) - \lambda X^{(2)}(t)\}\d t+\si\d W(t),\end{cases},$$
were $\lambda>0$ is a constant, $\si$ is an invertible $m\times m$-matrix, $W(t)$ is the $m$-dimensional Brownian motion, and $Z:\R^{2m}\times\scr P_2(\R^{2m})\to\R^m$ satisfies
$$|Z(x_1,\nu_1)-Z(x_2,\nu_2)|\le \aa_1|x_1^{(1)}-x_2^{(1)}|+\aa_2|x_1^{(2)}-x_2^{(2)}|+\aa_3\W_2(\nu_1,\nu_2)$$ for some constants $\aa_1,\aa_2,\aa_3\ge 0$ and all
  $ x_i=(x_i^{(1)}, x_i^{(2)})\in \R^{2m}, \nu_i\in \scr P_2(\R^{2m}), 1\le i\le 2.$
  If
\beq\label{*PW} 4\lambda> \inf_{s>0}\big\{2\aa_3 s + \aa_3s^{-1} + 2\aa_2+\ss{4(1+\aa_1)^2 + (2\aa_2+\aa_3s^{-1})^2}\big\},\end{equation}
then $\{L_t^\nu\}_{\nu\in \scr B_{r,R}} \in LDP(J)$ for all $r,R>0$.

Indeed,   $b(x,\nu):= (x^{(2)}-\lambda x^{(1)}, Z(x,\nu)-\lambda x^{(2)})$ satisfies
\beg{align*} & 2 \<b(x_1,\nu_1)-b(x_2,\nu_2), x_1-x_2\> \\
&\le -2\lambda |x_1^{(1)}-x_2^{(1)}|^2 -2(\lambda-\aa_2) |x_1^{(2)}-x_2^{(2)}|^2 \\
&\qquad + 2 |x_1^{(2)}-x_2^{(2)}|\big\{(1+\aa_1)|x_1^{(1)}-x_2^{(1)}| +\aa_3\W_2(\nu_1,\nu_2)\big\}\\
&\le \aa_3 s \W_2(\nu_1,\nu_2)^2 -\{2\lambda-\dd(1+\aa_1)\}|x_1^{(1)}-x_2^{(1)}|^2\\
&\qquad-\{2\lambda-2\aa_2-\dd^{-1}(1+\aa_1)-\aa_3s^{-1}\}|x_1^{(2)}-x_2^{(2)}|^2,\ \ s,\dd>0\end{align*}
for all $x_1,x_2\in\R^{2m}$ and $\nu_1,\nu_2\in \scr P_2(\R^{2m}).$
Taking
$$\dd= \ff{2\aa_2+\aa_3s^{-1} +\ss{4(1+\aa_1)^2+ (2\aa_2+\aa_3r^{-1})^2}}{2(1+\aa_1)}$$
such that $\dd(1+\aa_1)= 2\aa_2+\dd^{-1}(1+\aa_1)+\aa_3 s^{-1},$
we see that {\bf $(H_6^1)$} holds for some $\kk_1>\kk_2$ provided $2\lambda-\dd(1+\aa_1)>\aa_3 s$ for some $s>0$, i.e. \eqref{*PW} implies {\bf $(H_6^1)$}.
Moreover,  it is easy to see that conditions in Remark \ref{Remark 2.1.}(1) hold, see also \cite{GW12,WZ13} for  Harnack inequalities and gradeint estimates on stochastic Hamiltonian systems which also imply the strong Feller and $\bar\mu$-irreducibility of $\bar P_t$.      Therefore, the claimed assertion follows from Theorem  \ref{T01}(1).

\subsection{Distribution dependent SPDE}

Consider the following distribution-dependent SPDE on a separable Hilbert space $\H$:
  \beq\label{EPdd} \d X(t)= \{AX(t)+b(X(t),\L_{X(t)})\}\d t +\si(\L_{X(t)})\d W(t),\end{equation}
  where $(A,\D(A))$ is a linear operator on $\H$, $b: \H\times\scr P_2(\H)\to\H$ and $\si:  \scr P_2(\H)\to\mathbb L(\tt\H;\H)$ are measurable, and
  $W(t)$ is the cylindrical Brwonian motion on $\tt\H$.
   We make the following assumption.

  \beg{enumerate} \item[$(H_6^2)$]
  ($-A, \mathscr{D}(A))$ is    self-adjoint  with
  discrete spectrum
 $0<\lambda_1\le \lambda_2\le \cdots $ counting  multiplicities such that
 $\sum_{i=1}^\infty  \lambda_i^{\gg-1}<\infty$ holds for some constant $\gg\in (0,1)$.

Moreover,  $b$ is Lipschitz  continuous on $\H\times\scr P_2(\H)$, $\si$ is bounded and there exist constants $\aa_1,\aa_2\ge 0$ with $\lambda_1>\aa_1+\aa_2$ such that
$$2\<x-y, b(x,\mu)-b(y,\nu)\>+\|\si(\mu)-\si(\nu)\|_{HS}^2 \le 2\aa_1 |x-y |^2  +2\aa_2 \W_2(\mu,\nu)^2$$ holds for all $x,y\in \H$ and $ \mu,\nu\in \scr P_2(\H).$
\end{enumerate}
According to Theorem \cite[Theorem 3.1]{RW20a}, assumption {\bf $(H_6^2)$} implies
  that for any $X(0)\in L^2(\OO\to\H,\F_0,\P)$, the equation \eqref{EPdd} has a unique mild solution $X(t)$. As before we denote by $X^\nu(t)$ the solution with initial distribution $\nu\in \scr P_2(\H)$, and write $P_t^*\nu=\L_{X^\nu(t)}$.
Moreover, by It\^o's formula  and  $\kk:=\lambda_1-(\aa_1+\aa_2)>0$, it is easy to see that $P_t^*$ has a unique invariant probability measure $\bar\mu\in \scr P_2(\H)$ and
\beq\label{EXPP} \W_2(P_t^*\nu,\bar\mu)\le \e^{-\kk t} \W_2(\nu,\bar\mu),\ \ t\ge 0.\end{equation}
Consider the reference SPDE
$$\d \bar X(t)= \{A\bar X(t)+b(\bar X(t),\bar\mu)\}\d t +\si(\bar\mu)\d W(t),$$ which is again well-posed
for any initial value $\bar X(0)\in L^2(\OO\to\H,\F_0,\P)$. Let $J$ be the Donsker-Varadhan level 2 entropy function for the Markov process $\bar X(t)$, see \cite[Section 3]{RW20a}. For any $r,R>0$ let
$$\scr B_{r,R}:=\big\{\nu\in \scr P(\H): \nu(\e^{|\cdot|^r})\le R\big\}.$$

\beg{thm}[\cite{RW20a}] \label{T02} Assume $(H_6^2)$. If there exist constants $\vv\in (0,1)$ and $c>0$ such that
\beq\label{CCO}  \<(-A)^{\gg-1}x, b(x,\mu)\> \le c+ \vv|(-A)^{\ff \gg 2}x|^2,\ \  x\in\D((-A)^{\ff \gg 2}),  \end{equation}
then  $\{L^\nu_t\}_{\nu\in \scr B_{r,R}}\in LDP_u(J)$   for all $r,R>0$. If moreover $\bar P_t$ is strong Feller and $\bar\mu$-irreducible for some $t>0$, then $\{L^\nu_t\}_{\nu\in \scr B_{r,R}}\in LDP(J)$   for all $r,R>0$.
\end{thm}

Assumption $(H_6^2)$ is standard to imply  the well-posedness of \eqref{EPdd} and the exponential convergence of $P_t^*$ in $\W_2$.
 Condition \eqref{CCO} is implied by
\beq\label{GG} |(-A)^{\ff \gg 2 -1}b(x,\mu)|\le \vv' |(-A)^{\ff \gg 2}x|+c',\ \ x\in\D((-A)^{\ff \gg 2})\end{equation}  for some constants $\vv'\in (0,1)$ and $c'>0$. In particular,   \eqref{CCO} holds if $|b(x,\mu)|\le c_1+c_2|x|$ for some constants $c_1>0$ and $c_2\in (0, \lambda_1).$

  \subsection{Path-distribution dependent SPDE with additive noise}

Let $\tt\H=\H$ and $\si\in \mathbb L(\H).$     Then \eqref{*1}  becomes
\beq\label{**0} \d X(t)= \big\{AX(t)+ b(X_t,\L_{X_t})\big\}\d t + \si \d W(t).\end{equation}
Below we consider this equation with  $\si$ being invertible and  non-invertible respectively.

\subsubsection{Invertible $\si$}
Since $\si$ is constant, we are able to establish LDP for $b(\xi,\cdot)$ being Lipshcitz continuous in $\W_p$ for some $p\ge 1$ rather than just for $p=2$ as in the last two results.

\begin{enumerate}
\item[$(H_6^3)$] $\si\in\mathbb L(\H)$ is constant and $(A,\D(A))$ satisfies the corresponding condition in {\bf $(H_2)$}. Moreover, there exist constants $p\ge 1$ and  $\aa_1,\aa_2\ge 0$ such that
$$|b(\xi,\mu)-b(\eta,\nu)|\le \aa_1\|\xi-\eta\|_\infty +\aa_2 \W_p(\mu,\nu),\ \ \xi,\eta\in\C, \mu,\nu\in \scr P_p(\C).$$
 \end{enumerate}

Obviously,  $(H_6^3)$ implies assumption {\bf (A)} in \cite[Theorem 3.1]{RW20a},  so that for any $X_0^\nu\in L^p(\OO\to\C,\F_0,\P)$ with $\nu=\L_{X_0^\nu}$,
the equation \eqref{**0} has a unique mild segment solution $X_t^\nu$ with
$$\E\Big[\sup_{t\in [0,T]} \|X_t^\nu\|_\infty^p\Big]<\infty,\ \ T>0.$$
 Let $P^*_t\nu=\L_{X_t^\nu}$ for $t\ge 0$ and $\nu\in \scr P_p(\C)$.

When $P_t^*$ has a unique invariant probability measure $\bar\mu\in \scr P_p(\C)$,  we consider the reference  functional SPDE
\beq\label{**} \d\bar X(t)= \big\{A\bar X(t) + b(\bar X_t, \bar\mu)\big\}\d t + \si\d W(t).\end{equation}
By \cite[Theorem 3.1]{RW20a}, this reference equation is well-posed for any initial value in $L^p(\OO\to \C,\F_0,\P)$.
For any $\vv,R>0$, let
$$\scr I_{\vv,R}=\big\{\nu\in \scr P(\C):  \nu(\e^{\vv\|\cdot\|_\infty^2}) \le R\big\}.$$

 \begin{thm}[\cite{RW20a}] \label{TL1}
 Assume  $(H_6^3)$.    Let $\theta\in [0,\lambda_1]$ such that $$\kk_p:=
\theta-(\aa_1+\aa_2)\e^{p\theta r_0}=  \sup_{r\in [0,\lambda_1]} \big\{r- (\aa_1+\aa_2)\e^{prr_0}\big\}.$$
 \beg{enumerate} \item[$(1)$] For any $\nu_1,\nu_2\in \scr P_p(\C)$,
 \beq\label{ES1} \W_p(P_t^*\nu_1, P_t^*\nu_2)^p\le \e^{p\theta r_0-p\kk_p t} \W_p(\nu_1,\nu_2)^p,\ \ t\ge 0.\end{equation} In particular, if $\kk_p>0$, then   $P_t^*$ has a unique invariant probability measure $\bar\mu\in \scr P_p(\C)$   such that
\beq\label{ES2} \W_p(P_t^*\nu, \bar\mu)^p\le \e^{p\theta r_0-p\kk_p t} \W_p(\nu,\bar\mu)^p,\ \ t\ge 0, \nu\in \scr P_p(\C).\end{equation}
 \item[$(2)$] Let $\si$ be invertible. If $\kk_p>0$ and $\sup_{s\in (0,\lambda_1]}(s-\aa_1\e^{sr_0})>0$, then  $\{L_t^\nu\}_{\nu\in \scr I_{\vv,R}}\in LDP(J)$ for any $\vv,R>0$, where $J$ is the Donsker-Varadhan level 2 entropy function for the Markov process $\bar X_t$ on $\C$.
     \end{enumerate} \end{thm}

\paragraph{Example 6.3.} For a bounded domain $D\subset \R^d$, let $\H=L^2(D;\d x)$ and $A=-(-\DD)^\aa$, where $\DD$ is the Dirichlet  Laplacian on $D$ and
$\aa>\ff d 2$ is a constant. Let $\si=I$ be the identity operator on $\H$, and
$$b(\xi,\mu)= b_0(\mu)+ \aa_1 \int_{-r_0}^0 \xi(r)\Theta(\d r),\ \ (\xi,\mu)\in \C\times\scr P_1(\C),$$
where $\aa_1\ge 0$ is a constant,   $\Theta$ is a signed measure on  $[-r_0,0]$  with total variation $1$ (i.e. $|\Theta|([-r_0,0])=1$), and
$b_0$ satisfies
$$|b_0(\mu)-b_0(\nu)|\le \aa_2 \W_1(\mu,\nu),\ \ \mu,\nu\in\scr P_1(\C)$$ for some constant $\aa_2\ge 0$. Then $(H_6^3)$ holds for $p=1$, and as shown in he proof of Example 1.1 in \cite{BWY15} that $$\lambda_1\ge \lambda:=\ff{(d\pi^2)^{\aa}}{R(D)^{2\aa}},$$
where $R(D)$ is the diameter of $D$.
Therefore, all assertions in Theorem \ref{TL1} hold provided
$$\sup_{r\in (0,\lambda]} \{r-(\aa_1+\aa_2)\e^{rr_0}\}>0.$$ In particular, under this condition   $\{L_t^\nu\}_{\nu\in \scr I_{\vv,R}}\in LDP(J)$  for any $\vv,R>1.$

\subsubsection{Non-invertible $\si$}

Let $\H=\H_1\times\H_2$ for two separable Hilbert spaces $\H_1$ and
$\H_2$, and consider the following path-distribution dependent SPDE for $X(t)=(X^{(1)}(t),X^{(2)}(t))$ on $\H$:
\beq\label{E1} \beg{cases} \d X^{(1)}(t)= \{A_1X^{(1)}(t)+  BX^{(2)}(t)\}\d t,\\
 \d X^{(2)} (t)= \{A_2X^{(2)}(t)+ Z(X_t,\L_{X_t})\}\d t +\si\d W(t),
\end{cases}
\end{equation} where  $(A_i,\D(A_i))$ is a densely defined closed linear operator  on $\H_i$  generating a $C_0$-semigroup $\e^{t A_i}$  ($i=1,2$), $B\in \mathbb L(\H_2; \H_1)$,
$Z: \C\mapsto \H_2$ is measurable, $\si\in\mathbb L(\H_2)$, and  $W(t)$ is the cylindrical
Wiener process on $\H_2$.  Obviously,   \eqref{E1} can be reduced to  \eqref{**0} by taking $A={\rm diag}\{A_1,A_2\}$ and using ${\rm diag}\{0,\si\}$ replacing $\si$,  i.e.  \eqref{E1}  is a special case of  \eqref{**0} with non-invertible $\si$.

For any $\aa>0$ and  $p\ge 1$,   define $$ \W_{p,\aa}(\nu_1,\nu_2):=\inf_{\pi\in\C(\nu_1,\nu_2)}
\bigg(\int_{\C\times\C} \big(\aa\|\xi_1^{(1)}-\xi_2^{(1)}\|_\infty+ \|\xi_1^{(2)}-\xi_2^{(2)}\|_\infty\big)^p \pi(\d\xi_1,\d\xi_2)\bigg)^{\ff 1 p}.$$
 We assume

\beg{enumerate}\item[$(H_6^4)$] Let $p\ge 1$ and $\aa>0$.
  ($-A_2, \mathscr{D}(A_2))$  is   self-adjoint  with
  discrete spectrum
 $0<\lambda_1\le \lambda_2\le \cdots $ counting  multiplicities such that $\sum_{i=1}^\infty \lambda_i^{\gg-1}<\infty$ for some $\gg\in (0,1)$. Moreover, $ A_1 \le \dd-\lambda_1 $   for some constant  $ \dd\ge 0$; i.e., $\<A_1 x,x\>\le (\dd-\lambda_1)|x|^2$ holds for all $x\in \D(A_1)$.

Next,  there exist constants $K_1,K_2> 0$ such that
\begin{equation*}\begin{split}
&|Z(\xi_1,\nu_1)-Z(\xi_2,\nu_2)|\\
&\le
K_1\|\xi_1^{(1)}-\xi_2^{(1)}\|_\infty +K_2 \|\xi_1^{(2)}-\xi_2^{(2)}\|_\infty+K_3 \W_{p,\aa}(\nu_1,\nu_2),~~(\xi_i,\nu_i) \in\C\times\scr P_p(\C).
\end{split}\end{equation*}
Finally,  $\si$ is invertible on $\H_2$, and there exists $A_0\in \mathbb L(\H_1; \H_1)$ such that for any $t>0$, $B\e^{tA_2}= \e^{tA_1}\e^{tA_0} B$ holds   and
$$
Q_t:=\int_0^t\e^{sA_0}BB^*\e^{sA_0^*}\d s
$$
is invertible on $ \H_1$.
\end{enumerate}

By \cite[Theorem 3.2]{RW20a} for $\H_0=\H_2$ and ${\rm diag}\{0,\si\}$ replacing $\si$,    $(H_6^4)$ implies that for any $X_0\in L^p(\OO\to\C,\F_0,\P)$ equation \eqref{E1} has a unique mild segment solution. Let $P_t^*\nu=\L_{X_t}$ for $\L_{X_0}=\nu\in \scr P_p(\C).$

\begin{thm}[\cite{RW20a}] \label{T03} Assume $(H_6^4)$ for some constants $p\ge 1$ and $\aa>0$ satisfying
\beq\label{K1}\aa\le \aa':=\ff 1 {2\|B\|} \big\{\dd-K_2+\ss{(\dd-K_2)^2+4K_1\|B\|}\big\},\end{equation}
where $\|\cdot\|$ is the operator norm. If
\beq\label{K2}   \inf_{s\in (0,\lambda_1]}  s\e^{-  sr_0}>   K_2+\aa'\|B\| +  K_3,\end{equation}
then $P_t^*$ has a unique invariant probability measure $\bar\mu$ such that
\beq\label{K**} \W_p(P_t^*\nu,\bar\mu)^2\le  c_1\e^{-c_2 t} \W_p(\nu,\bar\mu),\ \ \nu\in \scr P_p(\C), t\ge 0 \end{equation}  holds for some constants $c_1,c_2>0$, and
  $\{L_t^\nu\}_{\nu\in \scr I_{\vv,R}}\in LDP(J)$  for any $\vv,R>1$, where $J$ is the Donsker-Varadhan level 2 entropy function for the associated reference equation for $\bar X(t)$.
\end{thm}

\paragraph{Example 6.4.} Consider the following equation for $X(t)= (X^{(1)}(t), X^{(2)}(t))$ on $\H=\H_0\times\H_0$ for a separable Hilbert space $\H_0$:
$$\beg{cases} \d X^{(1)}(t)= \{\aa_1 X^{(2)}(t)- \lambda_1 X^{(1)}(t)\}\d t\\
\d X^{(2)}(t)= \{Z(X_t,\L_{X_t}) - A   X^{(2)}(t)\}\d t+  \d W(t),\end{cases}$$
where $\aa_1\in\R\setminus\{0\}$, $W(t)$ is the cylindrical Brownian motion on $\H_0$,
$A$ is a self-adjoint operator on $\H_0$ with discrete spectrum such that all eigenvalues $0<\lambda_1\le\lambda_2\le\cdots $    counting multiplicities satisfy
$$\sum_{i=1}^\infty \lambda_i^{\gg-1}<\infty$$ for some $\gg\in (0,1),$ and    $Z$ satisfies
$$|Z(\xi_1,\nu_1)-Z(\xi_2,\nu_2)|\le \aa_2\|\xi_1-\xi_2\|_\infty+\aa_3 \W_2(\nu_1,\nu_2),\ \ (\xi_i,\nu_i)\in \C\times\scr P_2(\C), i=1,2.$$
Let $$\aa = \ff 1 {2\aa_1} \Big(\ss{\aa_2^2+4\aa_1\aa_2}-\aa_2\Big).$$ Then $P_t^*$ has a unique invariant probability measure $\bar\mu\in\scr P_2(\C)$, and
  $\{L_t^\nu\}_{\nu\in \scr I_{R,q}}\in LDP(J)$  for any $R,q>1$ if
\beq\label{ASS} \inf_{s\in [0,\lambda_1]} s\e^{-s r_0}> \aa_2+\aa_1\aa+\ff{\aa_3}{1\land\aa}.\end{equation}
  Indeed, it is easy to see that assumption {\bf $(H_6^4)$} holds for $p=2$, $\dd=0, \|B\|=\aa_1, K_1=K_2=\aa_2$ and $K_3=\ff{\aa_3}{1\land\aa}$. So, we have $\aa=\aa'$ and \eqref{ASS} is equivalent to \eqref{K2}. Then the desired assertion follows from
Theorem \ref{T03}.

\section{Comparison theorem}

The order preservation of stochastic processes is a  crucial property for one to compare a complicated process with simpler ones, and a result to ensure this property is called $``$comparison theorem" in the literature. There are two different type order preservations, one is in the distribution (weak) sense and the other is in the pathwise (strong) sense, where the latter implies the former.
  The weak order preservation has been investigated for  diffusion-jump Markov processes in \cite{CW,W} and references within,  as well as a class of super processes in \cite{W1}. There are also plentiful results on the strong order preservation,   see, for instance,
   \cite{BY, GD, IW, M, O, YMY} and references within for comparison theorems on forward/backward SDEs (stochastic differential equations), with jumps and/or with  memory. Recently,   sufficient and necessary conditions have been derived in   \cite{HW} for the order preservation of  SDEs with memory.

On the other hand, path-distribution dependent SDEs  have been investigated in \cite{HRW}, see also \cite{W18} and references within for distribution-dependent SDEs without memory. In this section, sufficient and necessary conditions of the order preservations for path-distribution dependent SDEs are presented.

Let $r_0\ge 0$ be a constant and   $d\ge 1$ be a natural number. The path space
${\bf C}=C([-r_0,0];\R^d)$ is Polish under the uniform norm $\|\cdot\|_\infty$.
For any continuous map $f: [-r_0,\infty)\to \R^d$  and
$t\ge 0$,  let  $f_t\in {\bf C}$ be such that $f_t(\theta)=f(\theta+t)$ for $\theta\in
[-r_0,0]$. We call $(f_t)_{t\ge 0}$ the segment of $(f(t))_{t\ge -r_0}.$
Next, let $\scr P({\bf C})$ be the set of probability measures on ${\bf C}$ equipped with the weak topology.
Finally, let $W(t)$ be an $m$-dimensional Brownian motion on a complete
filtration probability space $(\OO, \{\F_{t}\}_{t\ge 0},\P)$.

We consider the following Distribution-dependent SDEs with memory:
\beq\label{E1N} \beg{cases} \d X(t)= b(t,X_t,\L_{X_t})\,\d t+ \sigma(t,X_t,\L_{X_t})\,\d W(t),\\
\d \bar{X}(t)= \bar{b}(t,\bar{X}_t,\L_{\bar{X}_t})\,\d t+ \bar{\sigma}(t,\bar{X}_t,\L_{\bar{X}_t})\,\d W(t),\end{cases}\end{equation}
where
$$b,\bar{b}: [0,\infty)\times \C\times \scr P({\bf C})\to \R^d;\ \ \si,\bar{\sigma}: [0,\infty)\times \C\times \scr P({\bf C})\to \R^d\otimes\R^m$$
are   measurable.

For any $s\ge0$ and $\F_s$-measurable ${\bf C}$-valued random variables $\xi,\bar\xi$, a solution to (\ref{E1N}) for $t\ge s$ with $(X_s,\bar X_s)= (\xi,\bar\xi)$ is a  continuous adapted process $(X(t),\bar X(t))_{t\ge s}$  such that for all $t\ge s,$
\beg{equation*}\beg{split} &X(t) = \xi(0)+ \int_s^t  b(r,X_r,\L_{X_r})\d r+ \int_s^t \si(r, X_r,\L_{X_r})\d W(r ),\\
&\bar X(t) = \bar\xi(0)+ \int_s^t \bar b(r,\bar X_r,\L_{\bar{X}_r})\d r+ \int_s^t\bar \si(r,\bar X_r,\L_{\bar{X}_r})\d W(r) ,\end{split}\end{equation*} where  $(X_t, \bar X_t)_{t\ge s}$ is the segment process of $(X(t), \bar X(t))_{t\ge s-r_0}$ with
$(X_s,\bar X_s)= (\xi,\bar\xi)$.

Following the line of  \cite{HRW}, we consider the class of probability measures of finite second moment:
$$\scr P_2({\bf C})=\bigg\{\nu\in \scr P({\bf C}):\nu(\|\cdot\|_\infty^2):=\int_{{\bf C}}\|\xi\|_\infty^2\nu(\d\xi)<\infty\bigg\}.$$
It is a Polish space under the Wasserstein distance
$$\mathbb{W}_2(\mu_1,\mu_2):=\inf_{\pi\in \C(\mu_1,\mu_2)}\left(\int_{{\bf C} \times {\bf C}}\|\xi-\eta\|_\infty^2\pi(\d\xi,\d\eta)\right)^{\ff 1 2},\ \ \mu_1,\mu_2\in\scr P_2({\bf C}),$$
where $\C(\mu_1,\mu_2)$ is the set of all couplings for $\mu_1$ and $\mu_2.$

To investigate  the order preservation, we make   the following assumptions.

\beg{enumerate} \item[$(H_7^1)$] (Continuity)
There exists an increasing function $\aa: \R_+\to \R_+$     such that for any $ t\ge 0; \xi,\eta\in {\bf C}; \mu,\nu\in
\scr P_2({\bf C})$,
 \begin{align*}
&|b(t,\xi,\mu)- b(t,\eta,\nu)|^2+|\bar{b}(t,\xi,\mu)- \bar{b}(t,\eta,\nu)|^2+\|\si(t,\xi,\mu)- \si(t,\eta,\nu)\|_{HS}^2\\
&+\|\bar{\si}(t,\xi,\mu)- \bar{\si}(t,\eta,\nu)\|_{HS}^2\le \aa(t)\big(\|\xi-\eta\|_{\infty}^2+   \W_2(\mu,\nu)^2\big).
\end{align*}
\item[$(H_7^1)$] (Growth)  There exists an increasing function $K: \R_+\to\R_+   $ such that
 \begin{align*}
  |b(t,0,\dd_0)|^2+ |\bar b(t,0,\dd_0)|^2+\|\si(t,0,\dd_0)\|_{HS}^{2}+\|\bar\si(t,0,\dd_0)\|^{2}_{HS}
 \le K (t),\ \ t\ge 0,
 \end{align*}where $\dd_0$ is the Dirac measure at point $0\in {\bf C}$.
\end{enumerate}

It is easy to see that these two conditions imply assumptions $(H1)$-$(H3)$ in   \cite{HRW},  so that by \cite[Theorem 3.1]{HRW},  for any $s\ge 0$ and $\F_s$-measurable ${\bf C}$-valued random variables $\xi,\bar\xi$ with finite second moment,
 the equation  (\ref{E1N}) has a unique solution  $\{X(s,\xi;t), \bar X(s,\bar \xi;t)\}_{t\ge s}$    with $X_s=\xi$ and $\bar X_s=\bar\xi$.
Moreover, the segment process $\{X(s,\xi)_t, \bar X(s,\bar\xi)_t\}_{t\ge s}$ satisfies
\beq\label{BDD} \E \sup_{t\in [s,T]} \big(\|X(s,\xi)_t\|_\infty^2+ \|\bar X(s,\bar \xi)_t\|_\infty^2\big)<\infty,\ \ T\in [s,\infty).\end{equation}

To characterize the order-preservation for solutions of \eqref{E1N}, we introduce the partial-order on ${\bf C}.$ For $x=(x^1,\cdots, x^d)$ and $ y=(y^1,\cdots, y^d)\in\R^d$, we write $x\le y$ if $x^i\le y^i$ holds for all $1\le i\le d.$
Similarly,  for $\xi=(\xi^1,\cdots,\xi^d)$ and $\eta=(\eta^1,\cdots,\eta^d)\in {\bf C}$, we write $\xi\le \eta$ if $\xi^i(\theta)\le \eta^i(\theta)$ holds for all $\theta\in [-r_0,0]$ and $1\le i\le d.$ A function $f$ on ${\bf C}$ is called increasing if $f(\xi)\le f(\eta)$ for $\xi\le \eta$. Moreover, for any $\xi_1,\xi_2\in {\bf C}$, $\xi_1\land\xi_2\in {\bf C}$ is defined by $$(\xi_1\land\xi_2)^i:=\min\{\xi^i_1,\xi_2^i\}, \ \  1\le i\le d.$$
For two probability measures $\mu,\nu\in\scr P({\bf C})$, we   write $\mu\leq \nu$    if $\mu(f)\leq \nu(f)$  holds for any increasing function $f\in C_b({\bf C})$.  According to \cite[Theorem 5]{KKO},    $\mu\leq \nu$  if and only if there exists $\pi\in\C(\mu,\nu)$ such that
 $\pi(\{(\xi,\eta)\in {\bf C}^2: \xi\le \eta\})=1.$

\beg{defn} The stochastic differential system  $(\ref{E1N})$ is called order-preserving, if  for any $s\ge 0$ and $\xi, \bar\xi\in L^2(\OO\to {\bf C}, \F_{s},\P)$ with $\P(\xi\le\bar\xi)=1$,   $$\P\big(X(s,\xi;t)\le \bar X(s,\bar\xi;t),\ t\ge s\big)=1.$$   \end{defn}

We first present the following sufficient conditions for the order preservation, which reduce back to the corresponding ones in \cite{HW} when the system is distribution-independent.

\beg{thm}  Assume $(H_7^1)$ and $(H_7^2)$.   The system \eqref{E1N} is order-preserving provided the following  two conditions are satisfied:
\beg{enumerate}
\item[$(1)$] For any  $1\leq i\leq d$, $\mu,\nu\in \scr P_2({\bf C})$ with $\mu\leq\nu$, $\xi,\eta\in {\bf C}$ with $\xi\leq \eta$ and $\xi^i(0)=\eta^i(0)$,
    \begin{equation*}
b^i(t,\xi,\mu)\leq \bar{b}^i(t, \eta, \nu),\ \ \text{a.e.}\ t\ge 0.
    \end{equation*}
\item[$(2)$] For a.e.\ $t\ge 0$ it holds: $\si(t, \cdot,\cdot)= \bar\si(t,\cdot,\cdot)$    and $\sigma^{ij}(t,\xi,\mu)=\si^{ij}(t, \eta,\nu)$  for any $1\leq i\leq d$, $1\leq j\leq m$,  $\mu,\nu\in \scr P_2({\bf C})$ and $\xi,\eta\in {\bf C}$ with $\xi^{i}(0)=\eta^i(0)$.
\end{enumerate}
\end{thm}

Condition (2) means that for a.e. $t\ge 0$, $\si(t,\xi,\mu)=\bar\si(t,\xi,\mu)$   and the  dependence of $\si^{ij}(t,\xi,\mu)$ on $(\xi,\mu)$ is only via $\xi^i(0)$.

On the other hand, the next result shows that    these conditions are also necessary if all coefficients are continuous on $[0,\infty)\times{\bf C}\times\scr P_2({\bf C})$, so that \cite[Theorem 1.2]{HW} is covered when the system is distribution-independent.

\beg{thm}  Assume  $(H_7^1)$,  $(H_7^2)$      and that $\eqref{E1N}$ is order-preserving  for any complete filtered probability space $(\OO,\{\F_t\}_{t\ge 0},\P)$ and $m$-dimensional Brownian motion $W(t)$ thereon.  Then  for any $1\leq i\leq d$, $\mu,\nu\in \scr P_2({\bf C})$ with $\mu\le \nu$, and $\xi,\eta\in {\bf C}$ with $\xi\le \eta$ and $\xi^i(0)=\eta^i(0)$, the following assertions hold:
\beg{enumerate}
\item[$(1')$]  $b^i(t,\xi,\mu)\leq \bar{b}^i(t, \eta, \nu)$ if  $b^i$ and $\bar b^i$ are continuous  at points $(t,\xi,\mu)$ and $(t, \eta,\nu)$ respectively.
\item[$(2')$] For any $1\le j\le m$, $\sigma^{ij}(t,\xi,\mu)=\bar\si^{ij}(t, \eta,\nu)$ if $\si^{ij}$ and $\bar \si^{ij}$ are continuous   at points $(t,\xi,\mu)$ and $(t, \eta,\nu)$ respectively.
\end{enumerate}
Consequently, when $b, \bar b, \si$ and $\bar \si$ are continuous on $[0,\infty)\times {\bf C}\times \scr P_2({\bf C})$, conditions $(1)$ and $(2)$ hold.
\end{thm}

\end{document}